\documentclass[12pt,leqno]{article}
\usepackage[dvips]{epsfig}
\usepackage{graphics}
\usepackage{amssymb}

\newcommand\qed{\hfill$\sqcap\kern-7.5pt\hbox{$\sqcup$}$}

\newcommand{\NN}{\mathbb{N}}
\newcommand{\RR}{\mathbb{R}}
\newcommand{\R}{\mathbb{R}}
\newcommand{\N}{\mathbb{N}}
\newcommand{\Sp}{\mathbb{S}}

\newtheorem{theo}{Theorem}
\newtheorem{prop}[theo]{Proposition}
\newtheorem{lem}[theo]{Lemma}
\newtheorem{cor}[theo]{Corollary}
\newtheorem{rem}[theo]{Remark}

\newcommand{\beqn}{\begin{equation}}
\newcommand{\eeqn}{\end{equation}}
\newcommand{\bear}{\begin{eqnarray}}
\newcommand{\eear}{\end{eqnarray}}
\newcommand{\bean}{\begin{eqnarray*}}
\newcommand{\eean}{\end{eqnarray*}}

\newcommand{\Cc}{\mathcal{C}}

\newcommand{\EE}{\mathcal{E}}
\newcommand{\OO}{\mathcal{O}}

\newcommand{\Qq}{\mathcal{Q}}

\newcommand{\e}{{\varepsilon}}

\newcommand{\eps}{\varepsilon}
\newcommand{\wto}{\rightharpoonup}

\begin{document}

\title{Cooling process for inelastic Boltzmann equations for hard
spheres, Part II: Self-similar solutions and tail behavior}

\author{ S. {\sc Mischler}$^1$, C. {\sc Mouhot}$^2$ }

\footnotetext[1]{Ceremade, Universit\'e Paris IX-Dauphine,
place du M$^{al}$ DeLattre de Tasigny, 75016 Paris, France.}

\footnotetext[2]{UMPA, \'ENS Lyon, 46, alle d'Italie
69364 Lyon cedex 07, France.}

\maketitle

\begin{abstract}
We consider the spatially homogeneous Boltzmann
equation for inelastic hard spheres, in the framework of so-called
{\it constant normal restitution coefficients}. We prove
the existence of self-similar solutions, and we give pointwise
estimates on their tail. We also give general estimates
on the tail and the regularity of generic solutions.
In particular we prove Haff's law on the rate of decay of temperature,
as well as the algebraic decay of singularities.
The proofs are based on the regularity study of a rescaled problem,
with the help of the regularity properties of the gain part of
the Boltzmann collision integral, well-known in the elastic case,
and which are extended here in the context of granular gases.
\end{abstract}



\tableofcontents

\vspace{0.3cm}

\section{Introduction and main results}\label{sec:intro}

\setcounter{equation}{0}
\setcounter{theo}{0}

\subsection{The model}\label{subsec:model}

%{\bf 
%REVOIR 1: page 4 : rescaled variables
%REVOIR 3: page 28 : decay of singularities
%REVOIR 4: citer le "lemme g\'eom\'etrique"
%RELIRE: page 5 : motivation
%REMETTRE le dessin page 13
%}

We consider the asymptotic behavior of inelastic hard spheres described
by the spatially homogeneous Boltzmann equation with a
{\it constant normal restitution coefficient} (see~\cite{MMR1}).
More precisely,
the gas is described by the probability density of particles
$f(t,v) \ge 0$ with  velocity $v \in \RR^N$ ($N \ge 2$)
at time $t \ge 0$, which undergoes the evolution equation
     \bear  \label{eqB1}
     {\partial f \over \partial t}& = & Q(f,f) \quad\hbox{ in }\quad
     (0,+\infty)\times \RR^N,\\
     \label{eqB2}
     f(0) & = &  f_{\mbox{\scriptsize{in}}} \quad\hbox{ in }\quad  \RR^N.
     \eear
\smallskip

The bilinear collision operator $Q(f,f)$ models the interaction of
particles by
means of inelastic binary collisions (preserving mass and momentum
but dissipating kinetic energy).
Denoting by $e \in [0,1]$ the (constant) {\em normal restitution
coefficient}, when $e\not =0$ 
we define the collision operator in strong formulation as 
     \beqn \label{Qinel}
     Q(g,f)(v) := \int_{\RR^N \times \Sp^{N-1}} \left( \frac{~'f ~'g_*}{e^2}
     - f g_* \right) \, |u| \, b(\hat{u}\cdot \sigma) \, d\sigma \, dv_*,
                      %\nonumber  &=&
                      %Q^+ (g,f) - Q^-(g,f),
     \eeqn
where we use notations from~\cite{GPV**} (a dual formulation shall be given 
in (\ref{Qplusweak}) which includes the case $e=0$).
%and we have introduced the
%splitting between gain part $Q^+$ and loss part $Q^-$.
Here $u=v-v_*$ denotes the relative velocity, $\hat{u}$ stands for 
$u/|u|$, and $~'v,~'v_*$ denotes the possible pre-collisional
velocities leading to post-collisional velocities $v,v_*$. They
are defined by
\beqn
\label{primev}
'v = \frac{v+v_*}2 + \frac{'u}2, \hspace{0.5cm} 'v_* = \frac{v+v_*}2 - \frac{'u}2, 
\eeqn
with $~'u = (1-\beta) u + \beta |u| \sigma$ and $\beta = (e+1)/(2e)$ 
($ \beta \in [1,\infty)$ since $e \in (0,1]$). The elastic case corresponds 
to $e=1$. 
The function $b$ in (\ref{Qinel}) is (up to a multiplicative factor)
the differential collisional cross-section while $B =  |u| \, b(\hat{u} \cdot \sigma)$ 
represents the rate of collision of particles with pre-collisional velocities $v,v_* \in \R^N$ 
giving rise to particles with post-collisional velocities $v',v'_* \in \R^N$ defined by (\ref{vprime}). In the sequel we assume 
that there exists $b_0,b_1 \in (0,\infty)$ such that
     \beqn \label{cutoff}
     \forall \, x \in [-1,1], \quad b_0 \le b(x) \le b_1,
     \eeqn
and that 
  \begin{equation}\label{MM:hypmts}
  b \mbox{ is nondecreasing and convex on } (-1,1).
  \end{equation}
Note that the ``physical'' cross-section for hard spheres is given by
(see \cite{Cerci?,GPV**})
   $$
   b(x) = \mbox{cst} \, (1-x)^{-{N-3 \over 2}},
   $$
so that it fulfills hypothesis (\ref{cutoff}) and (\ref{MM:hypmts}) when $N=3$.
The Boltzmann equation~(\ref{eqB1}) is complemented with an initial
datum~(\ref{eqB2}) which satisfies (for some $k \ge 2$)
     \beqn  \label{initialcond}
     0 \le f_{\mbox{\scriptsize{in}}} \in L^1_k(\RR^N), \qquad
     \int_{\RR^N} f_{\mbox{\scriptsize{in}}} \, dv  = 1, \qquad
     \int_{\RR^N} f_{\mbox{\scriptsize{in}}} \, v \, dv  = 0.
     \eeqn
(see Subsection~\ref{MMR:subsec:not} for the notations of functional spaces). 
Notice that assuming the two last moment
conditions in~(\ref{initialcond}) is no loss of generality, since we
may always reduce to that case by a scaling and translation argument
(see \cite[section 1.5]{GPV**} for instance).
\smallskip

As explained in~\cite{MMR1}, the operator~(\ref{Qinel})
preserves mass and momentum:
     \beqn
     {d \over dt} \int_{\RR^N} f \,
      \left(
       \begin{array}{ll}
       1 \\
       v
       \end{array}
      \right)
     \, dv = 0,
     \eeqn
while kinetic energy is dissipated
     \beqn \label{eqdiffEE}
     {d \over dt} \EE(f(t,\cdot)) = - D(f(t,\cdot))
     \eeqn
where the energy $\EE$ and the dissipation functional $D$ are given by
     \[
     \EE (f) = \int_{\RR^N} f(v) \, |v|^2 \, dv,
\quad D (f):=  \tau \, \int_{\RR^N \times \RR^N} f \, f_* \, |u|^3 \, dv \, dv_*.
     \]
Here the inelasticity coefficient $\tau$ is defined by $ \tau  := m_b \, \left( \frac{1-e^2}{4} \right)$
and the angular momentum $m_b$ is defined by
\[
     m_b := \int_{\Sp^{N-1}} \left( \frac{1-(\hat{u}\cdot\sigma)}{2} \right) \, b(\hat{u}\cdot\sigma) \, d\sigma
        = |\Sp^{N-2}| \, \int_0 ^\pi b(\cos \theta) \, \sin^2 \theta/2 \,  \sin^{N-2} \theta \, d\theta
     \]
(in order to get the second formula, we have set $\cos \theta = \hat{u}\cdot\sigma$). 
\smallskip

The study of the Cauchy theory and the cooling process
of~(\ref{eqB1})-(\ref{eqB2}) was
done in~\cite{MMR1} (where more general models were considered). The
equation is well-posed for instance in $L^1 _2$: for $0 \le 
f_{\mbox{\scriptsize{in}}} \in L^1 _2$,
there is a unique solution in 
$C(\RR_+; L^1 _2) \cap L^1(\RR_+;  L^1 _3)$ (see again Subsection~\ref{MMR:subsec:not} 
for the notations of functional spaces).
This solution is defined for all times. It preserves mass, momentum and
has a decreasing kinetic energy. The cooling process does not
occur in finite time, but asymptotically in large time, {\it i.e.}, the kinetic energy
is strictly positive for all times and the solution satisfies
       \[  %\label{asymptTc}
       \qquad
       \EE(t) \to 0 \quad \hbox{and}  \quad f(t,\cdot) \ \wto \  \delta_{v=0}
       \,\,\hbox{ in }\,\, M^1(\RR^N)\hbox{-weak *}
       \quad \hbox{when} \quad t \to + \infty,
       \]
where $M^1(\RR^N)$ denotes the space of probability measures on $\RR^N$.
We refer to~\cite{MMR1} for the proofs of these results.

\subsection{Introduction of rescaled variables} \label{subsec:resc}

Let us introduce some rescaled variables, in order to study
more precisely the asymptotic behavior of the solution.
This usual rescaling can be found in~\cite{BGP**} and~\cite{EBphysrev}
for instance.  We search for a rescaled solution $g$ of the form
\beqn \label{fKgV}
f(t,v) = K(t) \, g(T(t),V(t) v ),
\eeqn
where $K,T,V$ are time scaling functions to be determined, such that 
$K(0)=V(0)=1$ and $T(0)=0$ (same initial datum).  
We choose the scaling functions $K,V$ such that they are compatible 
with {\it self-similar solutions}, that is when $g$ does not depend on time: there exists  a profile
function $G$ such that 
\beqn \label{fKGV}
f(t,v) = K(t) \,G(V(t) v ).
\eeqn
In this case, the conservation of mass 
$$
\mbox{cst} = \int_{\RR^N} f (t,v) \, dv = {K(t) \over V(t)^N} \, \int_{\RR^N} G(w) \, dw
$$
implies $K(t) = V(t)^N$. The evolution equation (\ref{eqB1})
satisfied by $f$ implies therefore  
\beqn\label{eqGV}
V'(t) \, \nabla_v \cdot \left( v \, G \right) = Q(G,G) 
\eeqn
by using the following homogeneity property: for any function $g$ such that 
the collision operator is well-defined, 
\beqn\label{homogeneiteQ}
\forall \, \lambda \in \RR^*, \quad 
Q(g(\lambda \cdot),g(\lambda \cdot))(v) = \lambda^{-(N+1)} \, Q(g,g)(\lambda v)
\eeqn
(which is obtained by a homothetic change of variable). Equation (\ref{eqGV}) then implies that 
$V'(t) = \mbox{cst} =: c_* >0$. 
When the rescaled solution $g$ does depend on time, its evolution 
equation is 
\beqn\label{eqgTV} T'(t) \, V(t) \, \partial_t g = Q(g,g) - c_* \, \nabla_v ( v \, g).
\eeqn
We then choose $T$ such this equation is as simple as possible: 
$T'(t) \, V(t) =1$. Hence we deduce the natural choice of the scaling functions
\beqn\label{Kc0c1}
K(t) = (1+c_* \, t)^N, \hspace{0.4cm} T(t) = \frac{1}{c_*} \, \ln \left(1
+ c_* \, t \right) \hspace{0.4cm} V(t) = (1 +c_* \, t)
\eeqn
for some constant $c_*>0$. It is obvious that changing $c_*$ in the equation (\ref{eqgTV})
only amounts to the multiplication of $g$ by a positive constant and 
the multiplication of $T'$ by a positive constant. In the sequel we fix without restriction $c_*=1$. 

Summarizing, thanks to the equation (\ref{eqGV}) and to the rescaled *
variables defined by (\ref{fKGV}), (\ref{Kc0c1}), for any {\it self-similar profile} $G$ solution to the stationary equation
  \beqn \label{eqrescG}
  Q(G,G) -  \nabla_v \cdot (v G) = 0
    \eeqn
we may associated a {\it self-similar solution} $F$ to the original equation~(\ref{eqB1})  by setting 
 \[ F(t,v) = (1+t)^N \, G( (1+t) v). \]
 Moreover, $G$ is obviously a {\it stationary solution} to the rescaled evolution equation 
   \beqn \label{eqresc}
    \frac{\partial g}{\partial t} = Q(g,g) -  \nabla_v \cdot (v g).
    \eeqn
which is the equation associated to (\ref{eqB1}) making the change of variables (\ref{fKgV}), (\ref{Kc0c1}) 
(with $c_*=1$). Roughly speaking the re-scaling (\ref{fKgV}), (\ref{Kc0c1}) 
adds an anti-drift to the original  equation (\ref{eqB1}).

More generally, for any solution
$g$ to the Boltzmann equation in self-similar variables (\ref{eqresc}),
we associate a solution $f$ to the evolution problem~(\ref{eqB1}),
defining $f$ by the relation
\beqn\label{ftog}
f(t,v) = (1+t)^N \, g(\ln(1 + t),(1+t) v).
\eeqn
Reciprocally, for any solution $f$ to the Boltzmann
equation~(\ref{eqB1}), we associate a solution $g$ to  the evolution
problem~(\ref{eqresc}), defining $g$ by the relation
\beqn\label{gtof}
g(t,v) = e^{-Nt} \, f(e^t -1, e^{-t} v).
\eeqn

Given an initial datum
$f_{\mbox{\scriptsize{in}}} = g_{\mbox{\scriptsize{in}}} \in L^1 _2$, 
we know from~\cite{MMR1} that
there exists a unique solution of~(\ref{eqB1}) in $C(\R_+,L^1 _2) \cap
L^1 (\R_+,L^1 _3)$. Therefore, thanks to the changes of variables 
(\ref{ftog}), (\ref{gtof}), we deduce that there exists a unique 
solution $g$ to~(\ref{eqresc}) in $C(\R_+,L^1 _2) \cap
L^1 _{\mbox{{\scriptsize loc}}}(\R_+,L^1 _3)$. Moreover we have the following relations 
between the moments of $f$ and $g$: 
\beqn \label{momentgtof}
\forall \, t \ge 0, \quad
\left\{
\begin{array}{l}
\|g(t,\cdot)  \, |\cdot|^k \|_{L^1} = e^{k \, t } \,  \|f(e^t-1,\cdot ) \, |\cdot|^k 
\|_{L^1} \vspace{0.3cm} \\
\|f(t,\cdot)  \, |\cdot|^k \|_{L^1} = (1+t)^{-k} \,  \|g(\ln(1+t),\cdot)  \, |\cdot|^k 
\|_{L^1}.
\end{array}
\right.
\eeqn

\subsection{Motivation}

The use of Boltzmann inelastic hard spheres-like models to describe dilute, 
rapid flows of granular media started with the seminal physics paper \cite{Haff}, and a huge physics litterature has 
developed in the last twenty years. The study of granular systems in such regime is motivated by 
their unexpected physical behavior (with the phenomena of {\em collapse} --or ``{\em cooling effect}''--
at the kinetic level and {\em clustering} at the hydrodynamical level), their use to derive hydrodynamical
equations for granular fluids, and their applications.
Granular gases are composed of  {\em macroscopic grains}, and not microscopic molecules 
like in rarefied  gas dynamics. The grains have only contact interactions, which motivates  
physical modelization by hard spheres with inelastic dissipative 
features built in the collision mechanism. The model of inelasticity with 
a constant normal restitution coefficient studied in the present paper is one of the simplest 
such model  (for a more elaborated model, see for instance the so-called {\it visco-elastic 
hard spheres model} in~\cite{BriPoe}, as well as \cite{MMR1}). 

From the physical and mathematical viewpoint, works on the inelastic Boltzmann models have been 
first restricted to the so-called {\em inelastic Maxwell molecules model}, which can be 
interpreted as an approximation where the collision rate is replaced by a 
mean value independent on the relative velocity. 
Existence, uniqueness of solutions and the rate of decay of the kinetic 
energy were obtained in~\cite{BCG00}  for the {\em inelastic Maxwell molecules model} 
with constant normal restitution coefficient.  
The {\em Maxwell molecules model} is important because of its analytic simplifications (with regards 
to the hard spheres model) allowing to use powerful Fourier transform tools. 
Polynomial tail behaviors of the self-similar profiles have been formally computed in \cite{EBJSP}. 
Convergence to self-similarity has been established in \cite{BC03,BoCeTo:03}. 
In the all, the {\em inelastic Maxwell model} is well understood now.
Similar results have also been obtained for some 
simplified non-linear friction models in \cite{CaglioVill,LiTo:04}. 

For the {\em inelastic hard spheres} model (with constant normal restitution coefficient), 
one easily sees from the discussion in Subsection~\ref{subsec:resc} 
that the kinetic energy of self-similar solutions 
to equation (\ref{eqB1}) (assuming their existence) behaves like 
\[ 
\EE(t) \mathop{\sim}_{t \to \infty} \frac{C}{t^2} \, .
\]
It is natural to expect a similar behavior for the rate of decay of the temperature for the generic 
solutions to equation (\ref{eqB1}). This conjecture was made twenty years ago in the pioneering paper~\cite{Haff}, 
and this rate of decay for the temperature is therefore known as {\em Haff's law}. Such a  law is a 
typical physical feature of inelastic hard spheres which does not hold for inelastic Maxwell molecules model 
(for which the temperature follows an exponential law).
Let us emphasize that in~\cite{BCG00} a {\em pseudo-Maxwell molecules model} was considered 
(multiplying the collision operator by some well-chosen scalar function of time), restoring  
Haff's law and still preserving the nice simplications of Maxwell molecules. However  
some important aspects of inelastic hard spheres such as the tail behavior are not preserved 
by this model.  

On the basis of the study of the particular case of Maxwell molecules, Ernst and Brito~\cite{EBphysrev} 
conjectured that self-similar solutions, when they exist, should attract any solution, in the sense of 
convergence of the rescaled solution. They also conjectured and formally computed 
some ``over-populated'' tail behaviors for the self-similar profile 
(namely decreasing slower than the Maxwellian) depending on the 
collision rate (see~\cite{EBJSP} for instance). 

Recently the works~\cite{GPV**,BGP**} laid the first steps for a
mathematical analysis of the inelastic hard spheres model 
(with constant normal restitution coefficient). In~\cite{GPV**}, the authors 
proved the existence of steady states and gave estimates showing 
the presence of over-populated tails for
{\em diffusively excited inelastic hard spheres}, that is when one provides 
an input of kinetic energy to the system preventing the collapse, 
which is modelized by some $\Delta_v f$ term added to equation~(\ref{eqB1}). 
In~\cite{BGP**} {\em a priori} integral estimates on the tail of the steady state (assuming
its existence) were established  for the spatially homogeneous inelastic Boltzmann
equation with various additional terms, such as a diffusion, or an anti-drift
as in~(\ref{eqresc}). 

In the present paper, we prove, for spatially homogeneous inelastic 
hard spheres with constant normal restitution coefficient, 
the existence of smooth self-similar solutions,
and we improve the estimates on their tails of~\cite{BGP**} into 
pontwise ones. 
We also give a complete regularity study of the generic solutions in the
rescaled variables, as well as estimates on their tails. In particular, 
we give the first mathematical proof of Haff's law 
and we show the algebraic decay of singularities.

In a forthcoming work~\cite{MMIII}, we shall prove the uniqueness and 
the local stability of these self-similar solutions for a small inelasticity.

\subsection{Notation}\label{MMR:subsec:not}

Throughout the paper we shall use the notation
$\langle \cdot \rangle = \sqrt{1+|\cdot|^2}$.
We denote, for any $q \in \RR$, the Banach space
     \[ 
     L^1_q = \left\{f: \RR^N \mapsto \RR \hbox{ measurable} \, ; \; \;
     \| f \|_{L^1_q} := \int_{\RR^N} | f (v) | \, \langle v \rangle^q \, dv
     < + \infty \right\}.
     \]

More generally we define the weighted Lebesgue space $L^p _q (\R^N)$
($p \in [1,+\infty]$, $q \in \R$) by the norm
    \[ \| f \|_{L^p _q (\R^N)} = \left[ \int_{\R^N} |f (v)|^p \, \langle v
        \rangle^{pq} \, dv \right]^{1/p} \]
when $p < +\infty$ and
    \[ \| f \|_{L^\infty _q (\R^N)} = \mbox{supess} _{v \in \R^N} |f (v)| \, 
         \langle v \rangle^{q} \]
when $p = +\infty$ (where $\mbox{supess}$ denotes the essential supremum).

The weighted Sobolev space $W^{k,p} _q (\R^N)$ ($p \in [1,+\infty]$, $q \in \R$ and $k \in \N$)
is defined by the norm
    \[ \| f \|_{W^{k,p} _q (\R^N)}  =   \left[ \sum_{|s| \le k} \|\partial^s f\|_{L^p _q} ^p \right]^{1/p} \]
where $\partial^s$ denotes the partial derivative associated with the  
multi-index $s \in \NN^N$. In the particular case $p=2$ we denote 
$H^k _q=W^{k,2} _q$. Moreover this definition can be extended to $H^s _q$ 
for any $s \ge 0$ by using the Fourier transform. 

We also denote by $L^1 _{\mbox{{\scriptsize loc}}}(\Omega)$ 
the space of locally integrable functions on a given set $\Omega \subset \RR^N$, 
that is the space of measurable functions on $\Omega$ 
which are integrable on every compact subset of $\Omega$.    
%We also introduce the space of normalized probability measures on
%$\RR^N$, denoted by $M^1 (\RR^N)$.
Finally, for $h \in \R^N$, we define the translation operator $\tau_h$ by
     \[ \forall \ v \ \in \ \R^N, \ \ \ \tau_h f (v) = f(v-h), \]
  %and for a vector
%$x \in \R^N$, we shall denote $\hat{x}=x/|x|$.
and we shall denote by ``$C$'' various constants which do not depend on the
collision kernel $B$.

\subsection{Main results}

In this subsection we consider a normal restitution coefficient 
$e \in (0,1)$ (except in dimension $N=3$ where the case $e=0$ can be included, 
see the discussions in the proofs). 

First we state a result of existence of self-similar solutions.
    \begin{theo}\label{selfsim}
    For any mass $\rho>0$, there exists a self-similar
    profile $G$ with mass $\rho$ and momentum $0$:
        \[ 0 \le G \in L^1 _2, \ \ \ \ Q(G,G) = \nabla_v \cdot (v \, G), \ \ \ \ 
          \int_{\R^N} G \left( \begin{array}{ll} 1 \\ v \end{array} \right) \, dv 
          = \left( \begin{array}{ll} \rho \\ 0 \end{array} \right), \]
    which moreover can be built in such a way that $G$ is radially
    symmetric, $G \in C^\infty$ and
        \[ \forall \, v \in \RR^N, \quad a_1 e^{- a_2 |v|} \le G(v)
             \le A_1 e^{- A_2 |v|} \]
      for some explicit constants $a_1,a_2,A_1,A_2 >0$.
% ici je propose juste de mettre C^\infty, le reste est
% implique par la borne exponentielle
   \end{theo}

Second we establish that Haff's law holds. 

   \begin{theo}\label{Haff} 
   For any $p \in (1,+\infty)$, $\rho >0$, and some initial datum 
   $f_{\mbox{\scriptsize{{\em in}}}}$ such that 
       \[ 0 \le f_{\mbox{\scriptsize{{\em in}}}} \in L^1 _2  \cap L^p, \ \ \
          \int_{\R^N} f_{\mbox{\scriptsize{{\em in}}}} \left( 
          \begin{array}{ll} 1 \\ v \end{array}\right) \, dv 
          = \left( \begin{array}{ll} \rho \\ 0 \end{array} \right), \]
   the associated solution of the Boltzmann equation~(\ref{eqB1},\ref{eqB2}) 
   in $C(\RR_+; L^1 _2) \cap L^1 (\RR_+;  L^1 _3)$ 
   satisfies Haff's law in the sense:
         \beqn \label{hafflaw}
         \forall \, t \ge 0, \hspace{0.3cm} \frac{m}{(1+t)^{2}} \le \EE(t) \le \frac{M}{(1+t)^{2}}
         \eeqn
       for some explicit constants $m,M>0$ depending on the collision 
       kernel and the mass, kinetic energy and $L^p$ norm of $f_{\mbox{\scriptsize{{\em in}}}}$.
%depending on the collision
%kernel and
%      the mass, energy and $L^p$ norm of $f_{\mbox{\scriptsize{in}}}$, 
%and $\tau$.
    \end{theo}

Third we give a more precise and general result on the regularity and asymptotic behavior 
of the solutions $g$ to the rescaled equation~(\ref{eqresc}) (note that the first point 
of Theorem~\ref{asymp} implies in particular the preceding theorem). 

\begin{theo}\label{asymp}
For any $p \in (1,\infty)$, $\rho>0$, and some initial datum $g_{\mbox{\scriptsize{{\em in}}}}$ such that 
      \[  
      0 \le g_{\mbox{\scriptsize{{\em in}}}} \in L^1_2  \cap L^p, \ \ \  
      \int_{\R^N} g_{\mbox{\scriptsize{{\em in}}}} \left(  \begin{array}{ll} 1 \\ v \end{array} \right) \, dv 
       =  \left( \begin{array}{ll} \rho \\ 0 \end{array} \right), 
      \]
the unique solution $g$ in  $C(\RR_+; L^1 _2) \cap L^1 _{\mbox{{\scriptsize {\em loc}}}} (\RR_+;  L^1 _3)$ 
of~(\ref{eqresc}) with initial datum $g_{\mbox{\scriptsize{{\em in}}}}$  satisfies:
       \begin{itemize}
       \item[(i)] It remains bounded in $L^p$ for all times, with uniform bound 
       as $t$ goes to infinity. Similarly for any $q \ge 0$, if 
       $g_{\mbox{\scriptsize{{\em in}}}} \in L^p _q$, then the solution remains 
       bounded in $L^p _q$ for all times with uniform bound as $t$ goes to infinity.   
       As for the Sobolev norms, for any $s,q \ge 0$, there is $w >0$ such that 
       if $g_{\mbox{\scriptsize{{\em in}}}} \in H^s _{q+w}$, then the solution 
       remains bounded in $H^s _q$ for all times, with uniform bound as $t$ goes 
       to infinity.   
       \item[(ii)] For any arbitrarily large $s,q \ge 0$, there exists $ \lambda>0$ and 
       some decomposition $g = g^S+g^R$ of the solution $g$ such that $g_S \ge 0$ and 
         \[
          \sup _{t \ge 0}   \left\|g_t ^S \right\|_{H^s_q} < +\infty, \qquad
         \left\|g_t ^R \right\|_{L^1_2}
         = O\left(e^{-\lambda t}\right) .
         \]
       \item[(iii)] Concerning the tail behavior, we have the following lower  and upper bounds: 
       there are some explicit constants $a_1,a_2 >0$ such that 
         \[ 
         \forall \, v \in \RR^N, \hspace{0.4cm} \liminf_{t \to \infty} g(t,v) \ge a_1 \, e^{- a_2 |v|}, 
         \]
       and, for any $\tau >0$ and $s \in [0,1/2)$, there are some explicit constants
       $A_1,A_2 >0$ such that (appearance of exponential moments)
         \[ 
         \forall \, t \ge \tau, \hspace{0.4cm} \int_{\RR^N} g(t,v) \, e^{- A_1 |v|^s} \, dv \le A_2. 
         \]
       Moreover similar integral upper bounds with $s\in [1/2,1]$ are uniformly propagated in time if they 
       are satisfied for the initial datum. 
          \end{itemize}
     All the constants in this theorem can be computed in terms of the mass, kinetic energy and 
     the different norms assumed on $g_{\mbox{\scriptsize{{\em in}}}}$, and the parameters.
      \end{theo}
\begin{rem}
Note that point~(ii) of Theorem~\ref{asymp} implies, by coming back to the 
original variables, that for an initial datum 
$0 \le f_{\mbox{\scriptsize{{\em in}}}} \in L^1 _2 \cap L^p$, $p \in (1,\infty)$,  
the unique associated solution of~(\ref{eqB1}) in 
$C(\RR_+; L^1 _2) \cap L^1 (\RR_+;  L^1 _3)$ 
satisfies a similar decomposition as above, but where the remaining part decreases 
with polynomial and not exponential rate. Hence we have shown that the amplitude of the 
singularities decreases algebraically in the original variables. 

More precisely it is likely that, in the original variables,  
singularities far from $0$ decrease in fact exponentially fast, whereas those 
close to $0$ cannot decrease faster than polynomially (due to the fact that 
the damping effect of the loss part of the collision operator degenerates at this point). 
%Moreover the Duhamel 
%representation of the solution in rescaled variables suggests that singularities 
%cannot decrease faster than exponentially in time for solutions of~(\ref{eqresc}), 
%and therefore singularities cannot decrease faster than algebraically for solutions 
%of~(\ref{eqB1}). 
\end{rem}

\subsection{Method of proof}

The main tool in this paper is the regularity theory of the collision operator: 
we show that its gain part satisfies similar regularity properties 
as in the elastic case \cite{LiQ,WeQ,BD98,Lu98,MV**}.
Following the study in the elastic case in~\cite{MV**}, we deduce uniform 
propagation of Lebesgue norms for the solutions in the rescaled variables (\ref{eqresc}). 

A first consequence is that the temperature in the rescaled variables is 
uniformly bounded from below by some positive number as soon as the initial datum 
satisfies some $L^p$ bound. Translating this estimate in the original variables, it proves Haff's law.

A second consequence is the existence of self-similar profiles (or steady states) 
for the rescaled equation~(\ref{eqresc}), which provides existence of self-similar 
solutions for the original equation~(\ref{eqB1}). 
This existence result is proved by the use of
a consequence of Tykhonov's fixed point Theorem (see Theorem~\ref{GPV}), which
is an infinite dimensional (rough) version of Poincar\'e-Bendixon
Theorem on dynamical systems, see for instance~\cite[Th\'eor\`eme 7.4]{Balabane} 
or \cite{GPV**,EMRR}. It states that a 
semi-group on a Banach space ${\cal Y}$ with suitable continuity properties, and 
which stabilizes a  nonempty convex weakly compact subset, has a steady state 
inside this subset.
We apply it to the evolution semi-group of~(\ref{eqresc}) in the Banach
space ${\cal Y}=L^1_2$. The existence and continuity properties of 
the semi-group were proved in~\cite{MMR1} and the nonempty
convex weakly compact subset of nonnegative functions with fixed mass and momentum and  
bounded moments and $L^p$ norm, $p \in (1,+\infty)$ (for some bound big enough)
is stable along the flow thanks to the above uniform $L^p$ bounds  
in the rescaled variables. 

Still following  the study in the elastic case in~\cite{MV**}, we also deduce from 
the regularity properties of the collision operator the uniform propagation of Sobolev 
norms as well as the exponential decay of (the amplitude of) the singularities for the 
solutions in the rescaled variables (\ref{eqresc}). That straightforwardly implies 
the smoothness of self-similar profiles as well as the algebraic decay of (the amplitude of) 
singularities for solutions to the Cauchy problem in the original variables (\ref{eqB1})-(\ref{eqB2}).

Let us now turn to the study of the tail behavior. 
On the one hand, we prove lower pointwise estimates on the self-similar profiles 
by the mean of some elementary maximum principles arguments inspired from~\cite{GPV**}. 
On the other hand, we prove explicit lower pointwise estimates on generic solutions in self-similar variables using the
spreading effect of the evolution semi-group associated to (\ref{eqresc}) (in the 
spirit of~\cite{Ca32,PW97,Mo**}). Finally,  upper pointwise estimates on the self-similar profiles 
are obtained using moments estimates established in~\cite{BGP**} and elementary o.d.e. arguments.

\subsection{Weak and strong forms of the collision operator}

Under our assumptions on $b$, the function $\sigma \mapsto b(\hat{u} \cdot \sigma)$ 
is integrable on the sphere $\Sp^{N-1}$, and we can set without restriction
   \[ \int_{\Sp^{N-1}} b(\hat{u}\cdot\sigma) \, d\sigma
        = |\Sp^{N-2}| \, \int_0 ^\pi b(\cos \theta) \,  \sin^{N-2} \theta \, d\theta =1. \]
Thus we can write the classical splitting $Q=Q^+ - Q^-$
between gain part and loss part. The loss part $Q^-$ is
  \beqn\label{Q-gfgPhif}
  Q^-(g,f)(v) := \left( \int_{\RR^N} \, g(v_*) \, |v-v_*| \, dv_*\right) \, f(v) = (g * \Phi) f,
  \eeqn
where $\Phi$ denotes $\Phi(z) = |z|$. For any distribution $g$ satifying
the moment conditions $\int_{\R^N} g \, dv = 1$, $\int_{\R^N} g \, v \,
dv = 0$, we have (see for instance~\cite[Lemma~2.2]{MMR1})
\beqn\label{Lgv} (g * \Phi) \ge |v|. \eeqn
%Representations of the gain term

The gain part $Q^+$ is defined by
  \beqn \label{Qinelbis}
  Q^+(g,f)(v) := \int_{\RR^N \times \Sp^{N-1}} \frac{~'f ~'g_*}{e^2} \,  |u|
  \, b(\hat{u}\cdot \sigma) \, d\sigma \, dv_*.
  \eeqn
%with (from [GPV]) $b$ the cross-section associated to contact
%collisions (hard spheres)
%    $$
%    b(x) = C_b \left( 1- x \right)^{-{N-3 \over 2}}
%    $$
%and
%  $u=v-v_*$ denotes the
%relative velocity, and
%$~'v$ and $~'v_*$ denote the pre-collisional velocities, defined by
%    \[ 'v = \frac{v+v_*}2 + \frac{'u}2, \hspace{0.5cm} 'v = \frac{v+v_*}2
%- \frac{'u}2 \]
%with $~'u = (1-\beta) u + \beta |u| \sigma$ and $\beta = (e+1)/2e \in
%(1,\infty)$ since $e \in (0,1)$.

\medskip
In the sequel, we shall need two other representations. On the one hand from
\cite{Cerci?}, there holds: for any $\psi \in L^\infty _1$, $f,g \in L^1_2$
\beqn \label{Qplusweak}
\int_{\RR^N}  Q^+(g,f)(v) \, \psi(v) \, dv = \int_{\RR^N \times \RR^N \times 
\Sp^{N-1}} f \, g_* \,  |u| \, b(\hat{u}\cdot \sigma) \, \psi(v') 
\, d\sigma \, dv_* \,dv,
\eeqn
where $v'$ denotes the post-collisional velocity  defined by
\beqn\label{vprime}
v' = \frac{v+v_*}2 + \frac{u'}2, \hspace{0.5cm} u' = \frac{1-e}2 \, u
+\frac{1+e}2 \,  |u|  \, \sigma.
\eeqn

On the other hand, we shall establish
a Carleman type representation for granular gases:

%Indeed it is a straightforward consequence of the following proposition.

%\subsection{A Carleman representation for granular gases}
%Here we prove the following representation of the gain term, for a
%general
%kernel $B(v-v_*,\cos \theta)$.

  \begin{prop} \label{carlrep}
  Let $E^e_{v,\!\!\!~'\!v}$ be the
  hyperplan orthogonal to the vector $v- \!\!\! ~'v$ and passing through
  the point $\Omega(v, \!\! ~'v)$, defined by 
    \[ \Omega(v,\!\! ~'v) := v + (1-\beta^{-1}) \, (v-\!\!~'v)
     = \left(2- \beta^{-1}\right) v + \left(\beta^{-1} -1 \right) \! ~'v. \] 
 Then we have the  following representation of the gain term
 \beqn \label{carlQ} \qquad\,\,\,
    Q^+(g,f)(v) = {2^{N-1} \over  \beta^{N-1} e^{2}}  \int_{\!\!\! ~'v \in \RR^N}  \!\!
    \int_{'v_* \in E^e_{v,\!\!\!~'\!v}} \!\!\!\! \!\!
    |v-v_*|^{2-N} B \,  |'v-v|^{-1} \, ~'g_* ~'f \, d \! ~'v \,  dE( \! ~'v_*) .
 \eeqn
Recall that $\beta=(1+e)/(2e)$ and $ B := B(u,\sigma) = |u|\, b(\hat u \cdot \sigma)$.
  \end{prop}
  
\begin{rem}\label{defN3e0} 
Let us emphasize that in dimension $N=3$ the expression (\ref{carlQ}) 
of the gain term $Q^+$ makes sense even when $e=0$ (since $(e \beta)^2$ 
converges to $1/4$ and $\Omega(v, \! ~'v) = 2 \, v - \! \! ~'v$ for $e=0$), while the formula 
defining $Q^+$ in (\ref{Qinel}) seems to be singular when $e\to0$. 
Hence this Carleman representation allows to define a strong 
formulation of $Q^+$ for $e=0$, at least in the physical case of the dimension $N=3$. 
\end{rem}

\medskip\noindent
{\sl Proof of Proposition~\ref{carlrep}}.
We start from the basic identity
  \begin{equation}\label{eq:dirac}
  \frac{1}{2} \, \int_{\Sp^{N-1}} F(|u|\sigma - u) \, d\sigma
  = \frac{1}{|u|^{N-2}} \, \int_{\RR^N} \delta(2 \, x \cdot u + |x|^2) \, F(x) \, dx,
  \end{equation}
which can be verified easily by completing the square in the Dirac
function, taking the spherical
coordinate $x+u=r \, \sigma$ and performing the change of variable $r^2 =s$.
We have the following relations from (\ref{primev})
\beqn\label{primev2}
     \left\{
     \begin{array}{ll}
     'v = v + (\beta/2) \, \left( |u| \sigma - u \right) \vspace{0.2cm} \\
     'v_* = v_* - (\beta/2) \, \left( |u| \sigma - u \right)
     \end{array}
     \right.
\eeqn
and thus starting from the strong form of $Q^+$ we get
     \[ Q^+(g,f) = e^{-2} \, \int_{\RR^N \times \Sp^{N-1}} B \,
         f\Big(v+(\beta/2) \, \left( |u| \sigma - u \right)\Big)
             g\Big(v_* - (\beta/2) \, \left( |u| \sigma - u \right) \Big) \,
        dv_* \, d\sigma. \]
Applying~(\ref{eq:dirac}) yields
   \[ Q^+(g,f) = 2 \, e^{-2} \, \int_{\RR^N \times \RR^N} |u|^{2-N} \, B  \,
   \delta(2 \, x \cdot u + |x|^2) \,
             f(v+(\beta/2)x) g(v_* - (\beta/2)x) \, dv_* \, dx. \]
We do the change of variable $x \to \! ~'v = v+(\beta/2)x$ (with jacobian $(\beta/2)^N$). 
Then, keeping $'v$ fixed, we make the change of variable $v_* \to \! ~'v_*$ 
(with jacobian $1$ since $'v_*=v+v_*- \!\! ~ 'v$). This gives
  \[
  Q^+(g,f) = {2^{N+1} \over \beta^{N} e^{2}} \, \int_{\RR^N \times  \RR^N}
  |u|^{2-N} \, B \,  \delta(2 \, x \cdot u + |x|^2) \, f('v) g('v_*) \, d\! ~'v_* \, d \! ~'v.
  \]
Finally, keeping $'v$ fixed, we decompose orthogonally
the variable $'v_*$ as $v+ V_1\,n+ V_2$ with $V_1 = (\! ~'v_*-v)
\cdot n$, $n =('v-v)/|'v-v|$ and $V_2$ orthogonal to $(\! ~'v-v)$. 
Let us compute the Dirac function in the new coordinates. Since 
$x=(2/\beta) \, (\! ~'v-v)$ and $u=(v-v_*)$, 
  \bean 
  2 \, x \cdot u + |x|^2 &=& (4/\beta) \, (\!~'v-v) \cdot (v-v_*) + (4/\beta^2) \, |\!~'v-v|^2 \\
           &=& (4/\beta) \, \Big( (\beta^{-1} -1) \, |\!~'v-v|^2 + (\!~'v-v) \cdot (\!~'v-v_*) \Big). 
  \eean
From the momentum conservation $(\!~'v-v_*) = (v - \! ~'v_*)$ and the orthogonal decomposition above: 
  \[  2 \, x \cdot u + |x|^2 = 
    (4/\beta) \, \Big( (\beta^{-1} -1) \, |\!~'v-v|^2 - V_1 \, |\!~'v-v|  \Big).  \]           
Hence we obtain the following representation: 
     \bean
     &&  Q^+(g,f) = {2^{N+1} \over \beta^{N} e^{2}} \, \int_{\R \times \R^{N-1}  \times \RR^N}  
     |u|^{2-N} \, B \, \\
        &&  \qquad \delta\left(\frac{4|'v-v|}{\beta} \Big[
    (\beta^{-1}-1) |'v-v| -V_1 \Big] \right) \,  f('v) g(v+V_1\, n + V_2) 
     \, dV_1 \, dV_2 \, d \! ~'v.
     \eean
It remains to remove the Dirac mass: we use the obvious identity
  \[
  \int_{\RR} \delta\left(\frac{4|'v-v|}{\beta} \Big[
    (\beta^{-1}-1) |'v-v| -V_1 \Big] \right) \, F(V_1) \, dV_1
    = \frac{\beta}{4 \, |\!~'v-v|} \, F\big((\beta^{-1}-1) |'v-v| \big) 
  \]
to finally obtain representation (\ref{carlQ}). \qed

\medskip

The parametrization by the Carleman representation means that 
for $v$ and $'v$ fixed, the point $'v_*$ describes the hyperplan
orthogonal to $('v-v)$ and passing through the point $\Omega(v,~\!'v)$ on the line
determined by $v$ and $'v$.
Note that in the elastic case, $\Omega(v,~\!'v) = v$, whereas here
$\Omega(v,~\! 'v)$
is outside the segment $[v,~\!'v]$, which reflects the fact that for 
the pre-collisional
velocities, the modulus of the relative velocity is bigger
than $|v-v_*|$. In the limit case $e=0$, $\Omega(v,~\! 'v)=2 \, v - \! \!~'v$. 
%For instance if one would have computed a Carleman representation for
%an inelastic fictious gain term
%with post-collisional velocities instead, he would have find
%$\Omega(v,v')$ belonging
%to the segment $[v,v']$.

\begin{figure}[h]
  \begin{center}
  \includegraphics[scale=.9]{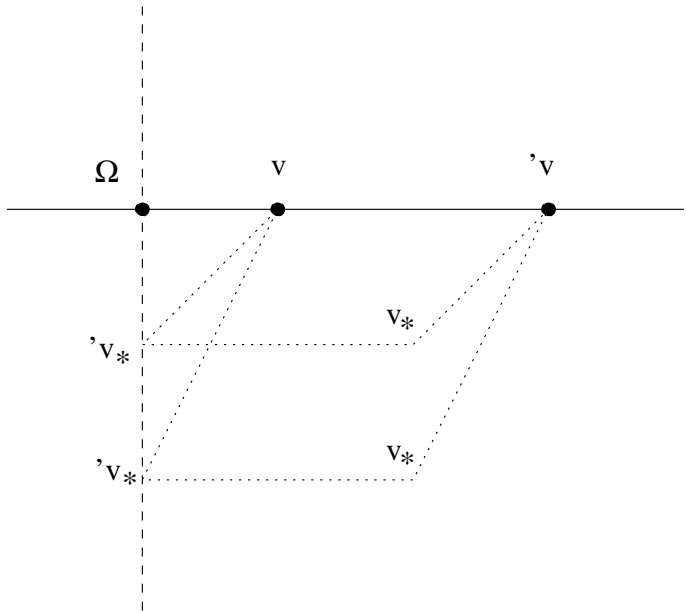}
  \caption{Carleman representation for granular gases}\label{carl}
  \end{center}
  \end{figure}

The geometrical picture (in a plane section) is summerized in
Figure~\ref{carl}.

From this proposition we immediately deduce the following representation, 
which is closer to the classical Carleman representation for 
the elastic Boltzmann collision operator.  From (\ref{primev2}) we deduce 
  \[ |\!~'v-v| = \frac{\beta}{\sqrt{2}} \, |u| \, \sqrt{1- (\hat u \cdot \sigma)}. \]
Hence we get 
  \[ %\label{QplusCarleman}
  Q^+(f,g)(v) = C_e \int_{\RR^N}  {~'f  \over |v-\!~'v|^{N-2}} \left\{
  \int_{E^e_{v,\!\!\!~'\!v}}~'g_*  \, \tilde b(\hat{u}\cdot \sigma) \,
  d \! ~'v_* \right\}\, d \! ~'v,
  \]
%where $E^e_{v,\!\!\!~'\!v}$ stands for the plane orthogonal to the 
%vector $'v-v$ and passing through the point
with 
  \[ C_e = \frac{ 2^{\frac{3N-5}2}}{\beta^{2N-4} e^2} \qquad  
%(de sorte que ce point $\Omega$ est exterieur au segment
%$[v,'v]$, attention au signe !!!! )
\mbox{ and } \qquad \tilde b(x) = (1-x)^{-(N-3)} \, b(x).  \]

\section{Regularity properties of the collision operator}\label{sec:Q}
\setcounter{equation}{0}
\setcounter{theo}{0}

In this section the final goal is to estimate quantities such as
     \[ \int_{\RR^N} Q(f,f) \, f^{p-1} \, dv \]
for $p>1$, {\it i.e.}, the action of the collision operator
on the evolution of the $L^p$ norm (to the power $p$) of the solution
along the flow. We shall use minoration estimates on
$Q^-$ deduced from (\ref{Q-gfgPhif})-(\ref{Lgv}), together with
convolution and regularity estimates on $Q^+$. The latters seem to be new
in the inelastic framework but they are an extension of
similar estimates in the elastic case $e=1$. 

The estimates on $Q^+$ can be splitted into tree groups. 
First the convolution-like estimates, which originated 
(in the elastic case) in the works of Gustafsson~\cite{G86,G88} 
(see also~\cite{MV**,DM,MMR1}) and were first extended to the inelastic 
case in~\cite{GPV**}. Second the regularization estimates in Sobolev spaces 
which originated in the elastic case 
in the works of Lions~\cite{LiQ}, Bouchut and Desvillettes~\cite{BD98},
Lu~\cite{Lu98} (see also~\cite{MV**} for some extensions). Third 
the non-concentration estimates in $L^1$, which originated in the elastic case 
in the work of Mischler and Wennberg~\cite{MW99}, and were extended by Abrahamsson~\cite{Abraham}. 
Let us also mention that regularity properties of $Q^+$ 
are reminiscent of the work of Grad on the linearized collision operator \cite{Grad63}. 
The main tool to extend 
the second group of estimates to the inelastic case shall be 
the Carleman representation for granular gases of
Proposition \ref{carlrep} (the third group of estimates can be extended 
with this tool as well, see~\cite{MMIII}). Before turning to the regularity
study of $Q^+$, we recall convolution-like estimates.

\subsection{Convolution-like estimates}\label{MM:subsec:conv}

In the elastic case $e=1$, convolution-like estimates for the gain part of
the collision operator were first proved in~\cite{G86,G88}.
This proof was simplified by a duality argument in~\cite{MV**},
where also a more precise statement was given.
These estimates were extended to the inelastic case, for
a constant normal restitution coefficient $e \in [0,1]$, in~\cite{GPV**} 
(in a form slightly less precise than in~\cite{MV**}). Also a result weaker in one 
aspect (less precise
for the treatment of the algebraic weight) but more general in another
(valid in any Orlicz spaces, and valid for more general collision
kernels) was proved in~\cite{MMR1}.
Here we only state the precise result we shall need, whose proof is
straightforward from the arguments in~\cite[Proof of Lemma~4.1]{GPV**}
and~\cite[Proof of
Theorem~2.1]{MV**}.

We make the following assumption on the cross-section:
no {\em frontal collision} should occur, {\it i.e.}, $b(\cos\theta)$
should vanish for $\theta$ close to $\pi$:
    \begin{equation}\label{eq:hypangconvol}
    \exists \, \theta_b>0 \ ; \quad
    \mbox{support} \, b \, ( \cos \theta) \subset \left\{ \theta \ /
    \ \ 0 \le \theta \le \pi - \theta_b \right\}. 
    \end{equation}
%This additional assumption will not be needed, however,
%for the quadratic estimates, {\it i.e.}, the estimates on $Q^+(f,f)$. Indeed,
%$Q^+(f,g) = \bar{Q}^+(g,f)$ if $\bar{Q}^+$ is the gain term 
%associated to the cross-section $\bar{b}(\cos \theta)
%= b(\cos (\pi - \theta))$. In particular, $b(\cos\theta)$ and
%$[b(\cos\theta) + b(\cos(\pi-\theta))] \, {\bf 1}_{\cos\theta\geq 0}$ define
%the same quadratic operator $Q^+$, and the latter
%satisfies~(\ref{eq:hypangconvol}) automatically (with $\theta_b =\pi/2$). 
%We note that $Q^+(g,f)$ and $Q^+(f,g)$ will not necessarily satisfy the same
%estimates, since assumption~(\ref{eq:hypangconvol}) is not symmetric.
To exchange the roles of $f$ and $g$, we introduce 
the symmetric assumption that no {\em grazing collision} should
occur, {\it i.e.}, 
    \beqn \label{eq:hypangconvol2}
    \exists \, \theta_b>0 \ ; \quad
    \mbox{support} \, b \, ( \cos \theta) \subset \left\{ \theta \ /
    \ \ \theta_b \le \theta \le \pi \right\}.
    \eeqn

Then we have (from the proofs of~\cite[Lemma~4.1 and Proposition~4.2]{GPV**}): 
    \begin{theo} \label{conv}
    Let $k,\eta \in \R$, $p\in [1, +\infty]$, and let $B = \Phi \, b$ be a
    collision kernel with $b$ satisfying the assumption~(\ref{eq:hypangconvol}).
    Then for any $e\in [0,1]$ the associated gain term satisfies the estimates
      \[ %\label{eq:conv}d
      \left\|Q^+(g,f)\right\|_{L^p _\eta} \le C_{k,\eta,p} (B)
      \left\|g\right\|_{L^1 _{|k+\eta|+|\eta|}}
      \left\|f\right\|_{L^p _{k+\eta}},
      \]
    with
      \[ C_{k,\eta,p} (B) = C \
          (\sin (\theta_b /2))^{\min(\eta, 0)-2/p'} \
          \left\|b\right\|_{L^1(\Sp^{N-1})}
          \left\|\Phi\right\|_{L^\infty _{-k}}.  \]
    If on the other hand assumption~(\ref{eq:hypangconvol}) is
    replaced by assumption~(\ref{eq:hypangconvol2}), then the same
    estimates hold with $Q^+(g,f)$ replaced by $Q^+(f,g)$.
    \end{theo}

\subsection{Lions Theorem for $Q^+$}

In this subsection we assume that the collision kernel $B= \Phi \, b$
satisfies
     \beqn \label{MM:hyplions}
     \Phi \in C^\infty _0 (\RR^N \setminus \{0\}), \hspace{0.4cm} b \in  C^\infty _0(-1,1).
     \eeqn
Then we have the
     \begin{theo}\label{regQ}
     Let $B$ be a collision kernel satisfying~(\ref{MM:hyplions}). Then for any 
     $e \in (0,1]$, the associated gain term $Q^+$ satisfies for all $s \in \R_+$ and $\eta \in \R_+$
       \[ \|Q^+(g,f)\|_{H^{s+(N-1)/2} _\eta} \le C(s,B) \, \|g\|_{H^s _\eta} \|f\|_{L^1 _{2\eta}} \]
       %\[ \|Q^+(g,f)\|_{H^{s+(N-1)/2} _\eta} \le C(s,B) \, \|g\|_{H^s _\eta} \|f\|_{L^1 _{2\eta}} \]
     for some explicit constants $C(s,B)>0$ depending only on $s$ and the collision kernel.
     \end{theo}

\begin{rem}
In dimension $N=3$ this theorem extends to the case $e=0$, with uniform bound 
$C(s,B)$ for $e \in [0,1]$. The only obstacle to the treatment of the case $e=0$ 
indeed is the constant $2^{N-1} \, \beta^{-(N-1)} \, e^{-2}$ in front of the Carleman 
representation, which may blow up as $e \to 0$. 
\end{rem} 

\medskip\noindent
{\sl Proof of Proposition~\ref{regQ}}.
We follow closely the proof of~\cite{MV**}, inspired from the works
of Lions~\cite{LiQ} and  Wennberg~\cite{WeQ}. Indeed the Carleman
representation proved above in Proposition~\ref{carlrep} allows essentially to reduce to the
study of the elastic case.

We assume first that $\eta=0$. We denote
     \[ {\cal B}(|\!~'v-\!~'v_*|,|\!~'v-v|) = \frac{B(|v-v_*|,\cos
          \theta)}{|v-v_*|^{N-2} |\!~'v-v|} \]
which belongs to $C^\infty_0((\R_+\setminus\{0\})^2)$ under
assumption~(\ref{MM:hyplions}).
We define the following (Radon transform type) functional: for $g$
smooth enough, $Tg$ is defined by
     \[ Tg(y) = \int_{\mu y + y^\bot} {\cal B}(z,y) \, g(z) \, dz \]
with $\mu = (2 - \beta^{-1})$.
Let us relate this functionnal with the Carleman representation~(\ref{carlQ}): we have 
     \bean && \int_{'v_* \in E^e_{v,\!\!\!~'\!v}} |v-v_*|^{2-N}  
     \, B \, |'v-v|^{-1} \, ~'g_* \, d\!~'v_*  \\  
     && \qquad = 
     \int_{'v_* \, \in \, \Omega(v,\!~'v) + (v-\!~'v)^\bot} \, {\cal B}(|\!~'v_*-\!~'v|,|v-\!~'v|) \, 
     \, ~'g_* \, d\!~'v_* \\
     && \qquad = \int_{z \, \in \, \Omega(v,\!~'v)-\!~'v + (v-\!~'v)^\bot} \, {\cal B}(|z|,|v-\!~'v|) \, 
     \, ~\tau_{-\!~'v}g(z) \, dz,   
     \eean
and since 
    \[  \Omega(v,\!~'v)-\!~'v = (2-\beta^{-1}) \, (v-\!~'v), \]
we deduce 
    \bean && \int_{'v_* \in E^e_{v,\!\!\!~'\!v}} |v-v_*|^{2-N}  
     \, B \, |'v-v|^{-1} \, ~'g_* \, d\!~'v_*  \\  
          && \qquad = \int_{\mu \, (v-\!~'v) + (v-\!~'v)^\bot} \, {\cal B}(|z|,|v-\!~'v|) \, 
     \, ~\tau_{-\!~'v}g(z) \, dz \\
     && \qquad = \left( \tau_{\!~'v} \circ T \circ \tau_{-\!~'v} \right) (g) (v).   
     \eean
Hence the representation~(\ref{carlQ}) writes 
     \[ Q^+(g,f)(v) = \frac{2^{N-1}}{\beta^{N-1} \, e^2} \, \int_{\RR^N} f(\! ~'v) \,
          \left( \tau_{\! ~'v} \circ T \circ \tau_{-\! ~'v} \right)(g)(v) \, d \! ~'v.\]
Thus if one has a bound on $T$ of the form
     \beqn \label{regT}
     \|Tg \|_{H^{s+(N-1)/2}} \le C_T \, \|g \|_{H^s}, \quad C_T >0, 
     \eeqn
then by using Fubini's and Jensen's theorems one gets
     \bean
     \|Q^+(g,f)\|^2 _{H^{s+(N-1)/2}} &\le& C \, \|f\|_{L^1} \, \int_{\RR^N} f(\!~'v) \,
          \left\| \left( \tau_{\!~'v} \circ T \circ \tau_{-\!~'v} \right)(g)\right\|^2 _{H^{s+(N-1)/2}} \, d\!~'v \\
          &\le& C \, \|f\|_{L^1} \, \int_{\RR^N} f(\!~'v) \,
          \left\| \left(T \circ \tau_{-\!~'v} \right)(g)\right\|^2 _{H^{s+(N-1)/2}} \, d\!~'v \\
         &\le& C \, C_T \, \|f\|_{L^1} \, \int_{\RR^N} f(\!~'v) \,
          \left\| \tau_{-\!~'v} g \right\|^2 _{H^s} \, d\!~'v  \\
         &\le& C \, C_T \, \left\| g \right\|^2 _{H^s} \,
            \|f\|_{L^1} \, \int_{\RR^N} f('v) \, d\!~'v
         \le C \, C_T \, \left\| g \right\|^2 _{H^s} \, \|f\|^2 _{L^1},
     \eean
which concludes the proof. Thus it remains to prove~(\ref{regT}).
But, up to an homothetic factor, $T$ is exactly the operator which was
studied in detail in~\cite{WeQ} and~\cite{MV**}. More precisely,
     \[ Tg(y) = \tilde{T} g (\mu y) \]
where $\tilde{T}$ is the Radon transform 
     \[ \tilde{T} g(y) = \int_{y + y^\bot} \tilde{{\cal B}} (z,y) \, g(z) \, dz, \]
which was introduced in the elastic case in \cite{WeQ}, associated 
with a kernel $\tilde {\cal B}$ related to our collision kernel by 
     \[ \tilde{{\cal B}}(z,y) =  {\cal B}(z, \mu^{-1} \, y). \]
It was proved in~\cite[Proof of Theorem~3.1]{MV**} that
     \[ \|\tilde{T} g \|_{H^{s+(N-1)/2}} \le C \, \|g \|_{H^s} \]
for an explicit bound $C$ depending on some weighted Sobolev 
norms on $\tilde{{\cal B}}$. Coming back to $T$, we
obtain~(\ref{regT}). This ends the proof when $\eta=0$.
The extension to $\eta >0$ is straightforward (and exactly similar 
to~\cite[Proof of Theorem~3.1]{MV**}). \qed

\medskip
As a Corollary we deduce from Theorem~\ref{regQ} the following estimate in
Lebesgue spaces by Sobolev embeddings (the proof is exactly similar
to \cite[Proof of Corollary~3.2]{MV**}).

    \begin{cor}\label{intQ}
    Let $B$ be a collision kernel satisfying~(\ref{MM:hyplions}).
    Then, for all $p \in (1,+\infty)$, $\eta \in \R$, we have
     \[
     \left\|Q^+(g,f)\right\|_{L^p _\eta} \le 
     C (p,\eta,B) \left\|g\right\|_{L^q _\eta} \left\|f\right\|_{L^1  _{2|\eta|}} 
     \]
     %\left\|Q^+(f,g)\right\|_{L^q _\eta} &\le&
     %C (p,\eta, B) \left\|g\right\|_{L^p _\eta} \left\|f\right\|_{L^1  _{2|\eta|}}
     %\eean
    where the constant $C (p,\eta,B) >0$ only depends on
    the collision kernel, $p$ and $\eta$, and $q<p$ is given by 
      \beqn \label{eq:defq}
      q= \left\{ \begin{array}{ll} \displaystyle 
      \frac{(2N-1)p}{N+(N-1)p} \quad \mbox{if $p \in (1,2N]$} \vspace{0.2cm} \\  \displaystyle 
      \frac{p}{N} \quad \mbox{if $p \in [2N,+\infty)$}. 
      \end{array} \right.
      \eeqn
     %\[
     %p      %\left\{
     %\begin{array}{cl} \displaystyle
     %\frac{q}{2-\frac{1}{N} + q \left(\frac{1}{N} -1 \right)}
     %& \mbox{if } q \in (1;2] \vspace{0.2cm} \\
     %qN & \mbox{if } q \in [2;+\infty).
     %\end{array}
     %\right.
     %\]
    \end{cor}

%\medskip\noindent
%{\sl Proof of Corollary~\ref{intQ}}.
%The proof is almost obvious. When $p=2$, it is a direct consequence of
%Theorem~\ref{regQ} with $s=0$, and the Sobolev injection
%$H^{\frac{N-1}{2}} _\eta \hookrightarrow L^{2N} _\eta$ (with a
%constant only depending on $N$). The general case follows by a
%Riesz-Thorin interpolation betweeen this estimate and
%Theorem~\ref{conv}. \qed

\subsection{Bouchut-Desvillettes-Lu Theorem on $Q^+$}

Now we turn to a slightly different regularity estimate on 
$Q^+$, which is a straightforward extension of the works~\cite{BD98,Lu98}
in the elastic case $e=1$. This class of estimate is weaker than
Lions's Theorem~\ref{regQ} since the Sobolev norm of $Q^+$ is controlled
by the square of the Sobolev norm of the solution with smaller order,
which does not allow to take advantage of the $L^1$ theory.
Nevertheless, it is more convenient in other aspects since it deals
directly with the physical collision kernel.

\begin{theo}\label{regQBDL} 
       Under the assumptions made on $B$ in Subsection~\ref{subsec:model}, for any $e \in [0,1]$ 
       the associated gain term $Q^+$ satisfies, 
       for all $s \in \R_+$ and $\eta \in \R_+$,
       \[ 
       \|Q^+(g,f)\|_{H^{s+(N-1)/2} _\eta} \le C(s,B) \, \left[ \|g\|_{H^s _{\eta+2}} \|f\|_{H^s _{\eta+2}}
       + \|g\|_{L^1 _{\eta+2}} \|f\|_{L^1 _{\eta+2}} \right] 
       \]
     for some explicit constant $C(s,B)>0$ depending only on $s$ and $B$.
     \end{theo}

\medskip\noindent
{\sl Proof of Theorem~\ref{regQBDL}}.
We follow closely the method in~\cite{BD98}. We write it for
$\eta=0$ but the general case is strictly similar.
%Let us consider the Fourier transform of $Q^+$.

Let us denote $F(v,v_*) = f(v) \, g(v_*) \, |v-v_*|$.
The same arguments as in~\cite{BD98} easily lead to
    \[ {\cal F}Q^+(\xi) = \int_{\Sp^{N-1}} \widehat F(\xi^+,\xi^-) \,
b(\hat{\xi} \cdot \sigma) \, d\sigma \]
where ${\cal F} Q^+$ denotes the Fourier transform of $Q^+$ according
to $v$, $\widehat F$ denotes the Fourier transform of $F$ according to $v,v_*$, and
    \[ %\left\{
       %\begin{array}{ll}
       \xi^+ = \frac{3-e}4 \, \xi + \frac{1+e}4 \, |\xi| \sigma, \qquad
       \xi^- = \frac{1+e}4 \, \xi - \frac{1+e}4 \, |\xi| \sigma.
       %\end{array}
       %\right.
    \]
Thus
    \[ |{\cal F}Q^+(\xi)|^2 \le \|b\|^2 _{L^2(\Sp^{N-1})} \,
       \left( \int_{\Sp^{N-1}} |\widehat F(\xi^+,\xi^-)|^2 \, d\sigma
       \right). \]
Let us consider frequencies $\xi$ such that $|\xi| \ge 1$. As
    \bean
    && \int_{\Sp^{N-1}} |\widehat F(\xi^+,\xi^-)|^2 \, d\sigma \\
       &=&  \int_{\Sp^{N-1}} \int_{|\xi|} ^{+\infty}
            -\frac{\partial}{\partial r}
    \left| \widehat F\left(\frac{3-e}4 \, \xi + \frac{1+e}4 \, r \sigma,
                   \frac{1+e}4 \, \xi - \frac{1+e}4 \, r \sigma\right) \right|^2 \, d\sigma \, dr \\
       &\le& C \, \int_{\Sp^{N-1}} \int_{|\xi|} ^{+\infty}
          \left| \widehat F\left(\frac{3-e}4 \, \xi + \frac{1+e}4 \, r \sigma,
                   \frac{1+e}4 \, \xi - \frac{1+e}4 \, r \sigma\right) \right| \times \\
       && \ \ \ \left|(\nabla_2-\nabla_1) \widehat F\left(\frac{3-e}4 \, \xi +
           \frac{1+e}4 \, r \sigma,
                   \frac{1+e}4 \, \xi - \frac{1+e}4 \, r \sigma\right) \right|
            \, d\sigma \, dr \\
       &\le& C \, \int_{|\zeta| \ge |\xi|}  \left| \widehat F\left(\frac{3-e}4 \,
        \xi + \frac{1+e}4 \, \zeta,
                   \frac{1+e}4 \, \xi - \frac{1+e}4 \, \zeta\right) \right| \, \times \\
       && \ \ \ \left|(\nabla_2-\nabla_1) \widehat F\left(\frac{3-e}4 \, \xi +\frac{1+e}4 \, \zeta,
                   \frac{1+e}4 \, \xi - \frac{1+e}4 \, \zeta\right) \right| \,
         \frac{d\zeta}{|\zeta|^{N-1}},
    \eean
where we have made the spherical change of variable $\zeta = r \sigma$, 
we deduce
    \bean
    && \int_{|\xi| \ge 1} |{\cal F}Q^+(\xi)|^2 \, |\xi|^{2s+(N-1)} \, d\xi \\
    && \le C \, \|b\|^2 _{L^2(\Sp^{N-1})} \,
        \int_{1 \le |\xi| \le |\zeta|} \left| \widehat F\left(\frac{3-e}4 \, \xi + \frac{1+e}4 \, \zeta,
                   \frac{1+e}4 \, \xi - \frac{1+e}4 \, \zeta \right) \right| \, \times \\
       && \ \left|(\nabla_2-\nabla_1) \widehat F\left(\frac{3-e}4 \, \xi + \frac{1+e}4 \, \zeta,
                   \frac{1+e}4 \, \xi - \frac{1+e}4 \, \zeta \right) \right| \,
                \frac{|\xi|^{2s+(N-1)}}{|\zeta|^{N-1}} \, d\xi \, d\zeta.
    \eean
Finally we make the change of variable
     \[ %\left\{
       %\begin{array}{ll}
       X = \frac{3-e}4 \, \xi + \frac{1+e}4 \, \zeta, \qquad
       Y = \frac{1+e}4 \, \xi - \frac{1+e}4 \, \zeta,
       %\end{array}
       %\right.
     \]
(whose Jacobian is uniformly bounded from above and below for $e \in [0,1]$) to obtain
    \bean \int_{|\xi| \ge 1} |{\cal F}Q^+(\xi)|^2 |\xi|^{2s+(N-1)} \, d\xi
         &\le& C \, \|b\|^2 _{L^2(\Sp^{N-1})} \,
        \int_{\RR^N \times \RR^N} \left| \widehat F(X,Y) \right| \times 
       \\ && \ \ \ \  \left|(\nabla_2-\nabla_1) \widehat F(X,Y) \right| \,
         \langle X \rangle^{2s} \, \langle Y \rangle^{2s}  \, dX \, dY \vspace{0.1cm} \\
      &\le& C \, \|b\|^2 _{L^2(\Sp^{N-1})} \, \|F\|_{H^s} \, \|(v-v_*)F\|_{H^s} \vspace{0.2cm} \\
      &\le& C \, \|b\|^2 _{L^2(\Sp^{N-1})} \, \|g\|_{H^s _2} ^2 \|f\|_{H^s _2} ^2.
    \eean
Then small frequencies are controlled thanks to the $L^1$ norms of $f$ and $g$,
which concludes the proof. \qed

\subsection{Estimates on the global collision operator in Lebesgue spaces}

We consider a collision kernel $B=\Phi \, b$ with
$\Phi(u)=|u|$ and $b$ integrable. We shall
make a splitting of $Q^+$ as in~\cite[Section~3.1]{MV**}.
We denote by ${\bf 1}_E$ the usual indicator function of the set $E$.

Let $\Theta: \R \rightarrow \R_+$ be an even $C^\infty$ function such
that $\mbox{support} \, \Theta \subset (-1,1)$, and $\int_\R \Theta \, dx= 1$. 
Let $\widetilde{\Theta}: \R^N \rightarrow \R_+$ be a radial
$C^\infty$ function such that $\mbox{support} \, \widetilde{\Theta} 
\subset B(0,1)$ and $\int_{\R^N} \widetilde{\Theta} \, dx = 1$. Introduce
the regularizing sequences
     \[
     \left\{
      \begin{array}{l}
      \Theta_m (x) = m \, \Theta(mx), \quad x\in \R, \vspace{0.2cm} \\
      \widetilde{\Theta}_n (x)       n^N \widetilde{\Theta}(n x), \quad x\in\R^N.
      \end{array}
     \right.
     \]
We use these mollifiers to split the collision kernel into a
smooth and a non-smooth part. As a convention, we shall use
subscripts $S$ for ``smooth'' and $R$ for ``remainder''. First, we
set
     \[ \Phi_{S,n} = \widetilde{\Theta}_n \ast
     \left( \Phi \ {\bf 1}_{{\cal A}_n} \right), \qquad
     \Phi_{R,n} = \Phi - \Phi_{S,n}, \]
where ${\cal A}_n$ stands for the annulus ${\cal A}_n = \left\{ x \in \R^N \ ; \ \frac{2}{n} \le |x| \le n \right\}$.
Similarly, we set
     \[ b_{S,m}(z) = \Theta_m \ast \left( b \ {\bf 1}_{{\cal I}_m} \right) (z), 
\qquad
       b_{R,m} = b - b_{S,m},\]
where ${\cal I}_m$ stands for the interval ${\cal I}_m = \left\{ x
\in \R \ ; \ -1+\frac{2}{m} \le |x| \le 1-\frac{2}{m} \right\}$
($b$ is understood as a function defined on $\R$ with compact
support in $[-1,1]$). Finally, we set
     \[ Q^+ = Q^+ _{S} + Q^+ _{R}, \]
where
     \[ %\label{eq:splitQ^+S}
     Q^+ _{S}(g,f) = e^{-2} \, \int _{\R^N \times \Sp^{N-1}}
     \Phi_{S,n}(|v-v_*|) \,
     b_{S,m} (\cos \theta) \, 'g_* \, 'f \, d\sigma \, dv_*
     \]
and
     \[ Q^+ _{R} = Q^+ _{RS} + Q^+ _{SR} + Q^+ _{RR} \]
with the obvious notation
     \[ \left\{
     \begin{array}{l} \displaystyle
     Q^+ _{RS} (g,f) = e^{-2} \, \int _{\R^N \times \Sp^{N-1}}
     \Phi_{R,n} \, b_{S,m} \, 'g_* \, 'f  \, dv_* \, d\sigma \vspace{0.3cm} \\ \displaystyle
     Q^+ _{SR} (g,f) = e^{-2} \, \int _{\R^N \times \Sp^{N-1}}
     \Phi_{S,n} \, b_{R,m} \, 'g_* \, 'f \, dv_* \, d\sigma \vspace{0.3cm} \\ \displaystyle
     Q^+ _{RR} (g,f) = e^{-2} \, \int _{\R^N \times \Sp^{N-1}}
     \Phi_{R,n} \, b_{R,m} \, 'g_* \, 'f \, dv_* \, d\sigma.
     \end{array} \right.
     \]

Now we follow the proof as in~\cite[Section~4.1]{MV**} since we have
the same functional inequalities in Sobolev and Lebesgue spaces, 
using also some ideas from~\cite{DM,MMR1} to simplify it. 
     \begin{prop}\label{gainLeb}
     Let us consider $e \in (0,1]$ and the associated gain term $Q^+$. 
     For any $\e >0$,  there exists $\theta\in (0,1)$, only
     depending on $N$ and $p$, and a constant $C_\e >0$, only depending on
     $N$, $p$, $B$ and $\e$ (and blowing up as $\e \to 0$), such that
       \[ \int_{\RR^N} Q^+(f,f) \, f^{p-1} \, dv \le
          C_\e \,  \|f\|_{L^1} ^{1+p\theta} \, \|f\|_{L^p} ^{p(1-\theta)}
                         + \e \, \|f\|_{L^1_2} \, \|f\|_{L^p _{1/p}} ^p. \]
     \end{prop}
\begin{rem} 
For general dimension $N \ge 2$, the estimates in this section 
are valid only for $e \in  (0,1]$ (the constants may blow up as $e \to 0$).  
However  when $N = 3$ they are uniform on $e \in (0,1]$ and then extends to the limit case $e=0$ 
(see also Remark~\ref{defN3e0}). 
\end{rem}

\medskip\noindent
{\sl Proof of Proposition~\ref{gainLeb}}.
Let us fix $\e >0$.
We split $Q^+$ as $Q^+ _{S} + Q^+ _{RS} + Q^+_{SR} +
Q^+ _{RR}$ and we estimate each term separately. 
%From the beginning we assume, without loss of generality, that the
%angular part $b(\cos\theta)$ of the collision kernel has its support
%included in $[0,\pi/2]$ (see the discussion on the symmetrization of $b$ 
%in Subsection~\ref{MM:subsec:conv}).
Remember that the truncation parameters $n$ (for the kinetic part)
and $m$ (for the angular part) are implicit in the decomposition
of $Q^+$.

By Corollary~\ref{intQ}, there exists a constant
$C(m,n)>0$, blowing up as $m$ or $n$ goes to infinity, such that
    \[ \left\|Q^+ _{S} (f,f)\right\| _{L^p}
    \leq C(m,n) \, \|f\|_{L^q} \|f\|_{L^1}, \] 
for some $q<p$ defined in~(\ref{eq:defq}). 
%    \beqn \label{eq:defq}
%    q= \left\{ \begin{array}{ll} \displaystyle 
%    \frac{(2N-1)p}{N+(N-1)p} \quad \mbox{if $p \in (1;2N]$} \vspace{0.2cm} \\  \displaystyle 
%    \frac{p}{N} \quad \mbox{if $p \in [2N;+\infty)$}. 
%    \end{array} \right.
%    \eeqn
%(the roles of $p$ and $q$ are exchanged here with respect to
%Corollary~\ref{intQ}).
Hence by H\"older's inequality,
    \[ \int_{\R^N} f^{p-1} Q^+ _S (f,f) \, dv
    \leq \left [ \int_{\R^N} f^p \, dv \right]^{\frac{p-1}{p}} \left
    [\int_{\R^N} (Q^+ _S)^p \, dv \right ]^{\frac1p} \]
    \[ \le \|f\|_{L^p}^{p-1} \big\|Q^+ _S(f,f)\big\|_{L^p} 
       \le   C(m,n) \, \|f\|_{L^q} \|f\|_{L^1} \, \|f\|_{L^p} ^{p-1}. \]

Next we fix a weight $\eta \ge -1$ and we estimate the $L^p _\eta$
norm of $Q^+ _{SR}(f,f)$ and $Q^+ _{RR}(f,f)$. We use that
$\|b_{R,m}\|_{L^1(\Sp^{N-1})}$ goes to $0$ as $m$ goes to infinity
(since $b$ is integrable on the sphere),
and we obtain, using Theorem~\ref{conv} with $k=1$ and splitting the angular integration 
between a part $\hat u \cdot \sigma \le 0$ with no grazing collision and a part 
$\hat u \cdot \sigma \ge 0$ with no frontal collision,
    \[ \big\|Q^+_{SR} (f,f), \ Q^+_{RR}(f,f) \big\|_{L^p_\eta} \leq \epsilon(m) \,
    \|f\|_{L^1_{|1+\eta|+|\eta|}} \|f\|_{L^p_{1+\eta}}, \]
for some $\epsilon(m)$ going to $0$ as $m$ goes to infinity.
%A similar estimate holds true for $\|Q^+_{SR}\|_{L^p_\eta}$. 
Since $1+\eta\geq 0$, we can write $|1+\eta|+|\eta| = 1+2\eta_+$, where $\eta_+= \max(\eta,0)$.
Hence by H\"older's inequality,
    \[ \int_{\R^N} f^{p-1} \Big( Q^+ _{SR} (f,f) + Q^+ _{RR}(f,f) \Big) \, dv
    = \int_{\R^N} \bigl (f \langle v\rangle^{1/p} \bigr )^{p-1}
    \frac{Q^+ _{SR} + Q^+ _{RR}}{\langle v\rangle^{\frac{1}{p'}}} \, dv \]
    \[ \leq \left [ \int_{\R^N}
    (f \langle v\rangle^{1/p})^p \, dv \right]^{\frac{p-1}{p}}
    \left[\int_{\R^N} \Big((Q^+ _{SR} + Q^+ _{RR} \Big) \langle v\rangle^{-1/p'} )^p \, dv \right]^{\frac1p} \]
    \[ \le  \|f\|_{L^p_{1/p}}^{p-1} 
          \Big( \|Q^+ _{SR} (f,f)\|_{L^p_{-1/p'}} + \|Q^+ _{RR} (f,f)\|_{L^p_{-1/p'}} \Big) 
       \le \epsilon(m) \, \|f\|_{L^1 _{1}} \, \|f\|_{L^p _{1/p}} ^p. \]

It remains to estimate the term corresponding to $Q^+_{RS}$. We have the trivial estimate 
  \[   \Phi_{R,n} \le C \, n^{-1} \, ( |v|^2 + |v_*|^2 ) \]
from which we deduce that 
  \bean  I &:=&  \int_{\R^N} f^{p-1} Q^+ _{RS}(f,f) \, dv \\ 
       &\le& C \, n^{-1} \, \int_{\R^N \times \R^N \times \Sp^N} b_{S,m} \, f \, f_* \, 
       (f')^{p-1} \, (|v|^2 + |v_*|^2) \, dv \, dv_* \, d\sigma := I_1 + I_2.
  \eean
Now we treat separately the two terms of the right-hand side: 
  \bean I_1 &=& \int_{\R^N \times \R^N \times \Sp^N} b_{S,m} \, (f \, |v|^2) \, f_* \, 
       (f')^{p-1} \, dv \, dv_* \, d\sigma \\
     &\le& C \, \int_{\R^N \times \R^N \times \Sp^N} b_{S,m} \, \left[ 
     f_* ^p + (f')^p \right] \, (f \, |v|^2) \, dv \, dv_* \, d\sigma := I_{1,1} + I_{1,2} 
   \eean
using Young's inequality $xy^{p-1} \le (1/p) \, x^p + ((p-1)/p) \, y^p$ on the product 
$f_* (f')^{p-1}$. 
The control of $I_{1,1}$ is immediate by integrating separately the angular variable: 
  \[ I_{1,1} \le C \, \| f \|_{L^1 _2} \, \|f \|_{L^p}. \]
For $I_{1,2}$,  as in the proof~\cite[Proposition 4.3]{GPV**} and using 
the notations of \cite[Lemma 4.4]{MMR1}, we make the change of 
variable $v_* \mapsto v' = \phi^*(v_*) = \phi^*_{e,v,\sigma}(v_*)$ keeping $v,\sigma$ 
fixed, which is a $C^\infty$-diffeomorphism from $\OO = \{ v_* \in \R^N, \hat u \cdot \sigma \not = 1 \}$ 
onto its image. Thanks to \cite[Lemma 4.4]{MMR1} and because $b_{S,m}$ 
has compact support in $(-1,1)$, its Jacobian $J^* = \hbox{det} (D\phi^*)$ satisfies
\beqn\label{JCm}
C(m)^{-1} \le J^*(v_*) \le C(m) \qquad \forall \, v_* \in \R^N, \ \hat u \cdot \sigma \in \hbox{support} \ b_{S,m},
\eeqn
for some constant $C(m) \in (0,\infty)$ which blows up when $m$ goes to $\infty$. 
Hence we straightforwardly deduce 
\beqn\label{I12m} \qquad
  I_{1,2} = C \int_{\R^N  \times \Sp^N}\! \int_{\phi^*(\OO )} \! b_{S,m} \, { (f')^p \over J ^*\circ \phi^{*-1}(v')}  
  \, f \, |v|^2 \, dv' \, dv \, d\sigma \le C(m) \, \| f \|_{L^1 _2} \, \|f \|_{L^p}.
\eeqn
The term $I_2$ is treated in a similar way: it is splitted 
as above using Young's inequality on the product $f (f')^{p-1}$. 
The  term $I_{2,1}$ involving $f^p$ is directly estimated as for the term $I_{1,1}$. 
For the term $I_{2,2}$ involving $(f')^p$ we proceed as for the term $I_{1,2}$. 
We make now the change of variable $v \mapsto v' = \phi_{e,v,\sigma}(v)$ 
keeping $v_*,\sigma$ fixed, where we use again the notations and results of 
\cite[Lemma 4.4]{MMR1}. Since its Jacobian $J$ also satisfies the bound (\ref{JCm}), 
we get the estimate (\ref{I12m}) for the term $I_{2,2}$. We  finally  deduce
  \[ 
   I \le \frac{C(m)}{n} \, \| f \|_{L^1 _2} \, \|f \|_{L^p}. 
   \]
Gathering the previous estimates we deduce 
    \bean 
    \int_{\R^N} f^{p-1} Q^+(f,f) \, dv &\leq& 
     C(m,n) \, \|f\|_{L^q} \, \|f\|_{L^1} \, \|f\|_{L^p}^{p-1} \\ 
    && \qquad 
    + \frac{C(m)}{n} \, \|f\|_{L^1_2} \, \|f\|_{L^p} 
    + \epsilon(m) \, \|f\|_{L^1_2} \, \|f\|_{L^p_{1/p}} ^p 
    \eean
where $q$ is defined by~(\ref{eq:defq}), $\epsilon(m)$ 
goes to $0$ as $m$ goes to infinity, and $C(m,n),C(m)>0$. Hence for any given 
$\e>0$, by first fixing $m$ big enough, then $n$ big enough, we get 
    \[ \int_{\R^N} f^{p-1} Q^+(f,f) \, dv \leq
     C_\e \, \|f\|_{L^q} \, \|f\|_{L^1} \, \|f\|_{L^p}^{p-1}
    + \e \, \|f\|_{L^1_2} \, \|f\|_{L^p_{1/p}} ^p \]
for some explicit constant $C_\e >0$. Combining this with
elementary interpolation, we deduce that there exists $\theta\in (0,1)$, only
depending on $N$ and $p$, and a constant $C_\e >0$, only depending on
$N$, $p$, $B$ and $\e$, such that
    \bean
    \int_{\R^N} f^{p-1} Q^+(f,f) \, dv & \leq C_\e
    \, \|f\|_{L^1} ^{1+p\theta} \, \|f\|_{L^p}^{1-p\theta} \, \|f\|_{L^p}^{p-1}
    + \e \, \|f \|_{L^1 _2} \, \|f\|_{L^p_{1/p}}^p \\
    & \leq C_\e \, \|f\|_{L^1} ^{1+p\theta} \, \|f\|_{L^p}^{p(1-\theta)}
    + \e \, \|f\|_{L^1_2} \, \|f\|_{L^p_{1/p}}^p.
    \eean
This concludes the proof.    \qed

\section{Regularity study in the rescaled variables}
\label{sec:self-sim}
\setcounter{equation}{0}
\setcounter{theo}{0}

In this section we show the uniform  propagation of Lebesgue and
Sobolev norms and the exponential decay of singularities for the
solutions of~(\ref{eqresc}).

\subsection{Uniform propagation of moments - Povzner Lemma}

Let us prove that the kinetic energy of $g$ remains uniformly bounded from above 
as $t$ goes to infinity. Using (\ref{eqresc}) and (\ref{eqdiffEE}), we get
    \[ {d \over dt}  \int_{\RR^N} g \, |v|^2 \, dv \le - \tau \, \int_{\RR^N \times \RR^N} 
       g \, g_* \, |u|^3 \, dv_*dv + 2 \,  \int_{\RR^N} g  \, |v|^2 \, dv. \]
On the one hand, from Jensen's inequality (see for instance \cite[Lemma~2.2]{MMR1}), there holds
    \[
         \int_{\RR^N} g_* \, |u|^3 \, dv_* \ge \rho \, |v|^3.
    \]
On the other hand,  H\"older's inequality yields
    \[  \int_{\RR^N} g \, |v|^2 \, dv \le \left(  \int_{\RR^N} g \, dv \right)^{1/3}
                                 \left(  \int_{\RR^N} g |v|^3 \, dv \right)^{2/3}, \]
which implies that
    \[  \int_{\RR^N} g \, |v|^3 \, dv \ge  \rho^{1/2} \, \left(
\int_{\RR^N} g |v|^2 \, dv \right)^{3/2}. \]
Thus
  \bean
  {d \over dt}  \int_{\RR^N} g \, |v|^2 \, dv \le
  - \tau \, \rho^{3/2} \, \left(  \int_{\RR^N}  g \, |v|^2 \, dv \right)^{3/2}
  + 2 \, \left( \int_{\RR^N}  g \, |v|^2 \, dv \right) \\
  \le  \tau \, \rho^{3/2} \, \left( \int_{\RR^N}  g \, |v|^2 \, dv \right)
  \left[ \frac{2}{ \tau \rho^{3/2}} - \left(  \int_{\RR^N}  g \, |v|^2 \, dv \right)^{1/2} \right],
  \eean
and by maximum principle we deduce
  \beqn\label{BorneY2}
  \sup_{t \ge 0}  \int_{\RR^N} g \, |v|^2 \, dv \le C_E             \max \left\{ \left(\frac{4}{\rho^3 \,  \tau^2}\right),
  \int_{\RR^N} g_{\mbox{\scriptsize{in}}} \, |v|^2 \, dv \right\}.
  \eeqn

The same argument, together with sharp versions of Povzner inequalities 
from~\cite{B97,BGP**}, yields uniform bounds and appearance on every moments of the solution, 
as well as appearance of some exponential moments (this last point 
was first noticed in~\cite{MMR1}), in a similar way as 
in~\cite[Proof of Proposition~3.2]{MMR1}. Indeed we prove the
      \begin{prop}\label{MM:Y3uniform}
      Let $g$ be a solution  in $C(\RR_+; L^1 _2) \cap 
      L^1 _{\mbox{{\scriptsize{\em loc}}}} (\RR_+;  L^1 _3)$ 
      to the rescaled Boltzman equation~(\ref{eqresc}) with $e \in [0,1)$, 
      with initial datum $g_{\mbox{\scriptsize{{\em in}}}}$.  
      Then it satisfies the following additional moment properties:
        \begin{enumerate}
        \item[(i)]
        For any $s \ge 2$, there is an explicit constant $C_s > 0$, depending only on $B$,
        $e$, and $g_{\mbox{\scriptsize{{\em in}}}}$, such that
          \[
          \sup_{t \in [0,\infty)} \| g(t,\cdot) \|_{L^1_s} \le
          \max\big\{ \| g_{\mbox{\scriptsize{{\em in}}}} \|_{L^1_s},C_s \big\}.
          \]
        \item[(ii)]
        If $g_{\mbox{\scriptsize{{\em in}}}} \, e^{r \, |v|^\eta} \in L^1(\R^N)$ for $r>0$ and $\eta \in (0,1]$,
        there exists $C_1, r' > 0$, depending only on $B$,
        $e$, and $g_{\mbox{\scriptsize{{\em in}}}}$, such that
          \[
          \sup_{t \in [0,\infty)} \int_{\R^N} g(t,v) \, e^{r' \, |v|^\eta} \, dv \le C_1.
          \]
       \item[(iii)] For any  $\eta \in (0,1/2)$ and $\tau >0$, there 
exists $a_\eta,C_\eta \in (0,\infty)$,
       depending only on $B$, $e$, $\tau$ and $g_{\mbox{\scriptsize{{\em in}}}}$,
       such that
          \[ %\label{expeta}
          \sup_{t \in [\tau,\infty)} \int_{\R^N} g(t,v) \, e^{a_\eta \,
          |v|^\eta} \, dv \le C_\eta.
          \]
         \end{enumerate}
Let us emphasize that the constant $C_s, a_\eta, C_\eta$ may depend on
$g_{\mbox{\scriptsize{{\em in}}}}$ only  through
its mass $\rho$ and its kinetic energy $\EE_{\mbox{\scriptsize{{\em 
in}}}}$.
\end{prop}

\smallskip\noindent
{\sl Proof of Proposition \ref{MM:Y3uniform}. } The proof is just a copy
with minor modifications of classical proofs. For the proof of (i)  we
refer for intance to \cite{MW99,CVhand,GPV**} and the references
therein. The proofs of (ii) and (iii) are variants of the proof
of~\cite[Proposition~3.2]{MMR1}, which itself follows closely the proof
of \cite[Theorem~3]{B97} extended to the inelastic case in
\cite{BGP**}. The starting point is the following differential
equation  on the moments
  \[ 
  {d \over dt} m_p = \int_{\R^N} Q(g,g) \, |v|^{2p} \, dv + p \, m_p \quad\hbox{with}\quad
  m_p := \int_{\R^N} g \, |v|^{2p} \, dv.
  \]
Using the same notation as in \cite[Proof of Proposition~3.2]{MMR1}, we 
introduce the new rescaled moment function
  \[ z_p := {m_p \over \Gamma(a\,p+1/2)}, \quad
     Z_p := \max_{k=1,\dots,k_p} \{ z_{k+1/2} \, z_{p-k}, \,  z_{k} \, z_{p-k+1/2} \}, \]
for some fixed $a \ge 2$,   and we obtain the  differential inequality
  \beqn \label{MM:mp10}
  {d z_p \over dt} \le A'  \,   p^{a/2-1/2} \, Z_p - A'' \,
  p^{a/2} \, z_p^{1+1/2p} + p \, z_p
  \eeqn
for any $p = 3/2,2, \dots$ and for some constants $A',A'' > 0$. Note 
that (\ref{MM:mp10}) is nothing but~\cite[equation (3.18)]{MMR1},
with an additional term $p \, z_p$ due to the additional term 
$-\nabla_v \cdot (v \, g)$ in
equation~(\ref{eqresc}).

\smallskip
On the one hand, we remark, by an induction argument, that taking $p_0
= p_0(a,A',A'')$ and $x_0 = x_0(a,A',A'')$ large enough, the sequence 
of functions
$z_p := x^p$ is a sequence of supersolution
of~(\ref{MM:mp10}) for any $x \ge  x_0$ and $p \ge p_0$. Let us emphasize
here that we have to take $a \ge 2$ ({\it i.e.}, $\eta \le 1$ in 
\cite[Proof of Proposition~3.2]{MMR1}) because of the additional term $p \, z_p$. 
On the other hand, choosing $x_1$
large enough, which  may depend on $p_0$, we have from (i)  that  the
sequence of functions $z_p := x^p$ is a
sequence of  supersolution of~(\ref{MM:mp10}) for any $x \ge x_1$ and
for $p \in \{ 0,1/2, ... , p_0 \}$. As a consequence,  we have proved
that there exists $x_2 := \max(x_0,x_1)$ such that  the set
  \beqn \label{MM:superS}
  \Cc_x := \left\{ z = (z_p); \quad z_p \le x^p \,\, \, \forall \, p
  \in {1 \over 2} \, \NN \right\}
  \eeqn
is invariant under the flow generated by the Boltzmann equation
for any $x \ge x_2$: if $g(t_1) \in \Cc_x$ then $g(t_2) \in \Cc_x$ for
any $t_2 \ge t_1$.
The end of the proof is exactly similar to that of~\cite[Proof of Proposition~3.2]{MMR1}.
\qed
\smallskip

The integral upper bound in point (iii) of Theorem~\ref{asymp} follows from 
point (iii) of Proposition~\ref{MM:Y3uniform}. 

\subsection{Stability in $L^1$}

The stability result~\cite[Proposition~3.3]{MMR1} translates
for~(\ref{eqresc}) into:
    \[ \|g-h\|_{L^1} + e^{-2T} \, \|(g-h)|v|^2\|_{L^1}
       \le e^{C(e^{2T}-1)} \left[ \|g_{\mbox{\scriptsize{in}}}-h_{\mbox{\scriptsize{in}}}\|_{L^1}
     +  \|(g_{\mbox{\scriptsize{in}}}-h_{\mbox{\scriptsize{in}}})|v|^2\|_{L^1} \right] \]
for any solutions $g$ and $h$ in $C(\R_+,L^1 _2) \cap L^\infty(\R_+,L^1 _3)$ with initial 
datum $0 \le g_{\mbox{\scriptsize{in}}}, h_{\mbox{\scriptsize{in}}} \in L^1 _3$. 
This shows that, in the Banach space $L^1 _2$, 
the evolution semi-group $S_t$ of~(\ref{eqresc}) satisfies: for any $t \ge 0$, 
$S_t$ is (strongly) continuous in any $L^1 _3$ bounded subset of $L^1 _2$. 
However we shall prove a more precise stability result, working directly on the rescaled
equation~(\ref{eqresc}).

    \begin{prop}\label{stab}
    Let $0 \le g_{\mbox{\scriptsize{{\em in}}}}, h_{\mbox{\scriptsize{{\em in}}}} \in L^1 _3$ and let 
    $g$ and $h$ be the two solutions of~(\ref{eqresc}) (in $C(\R_+,L^1 _2) \cap L^\infty(\R_+,L^1 _3)$) 
    with $e \in [0,1]$. Then there is $C_{\mbox{\scriptsize{{\em stab}}}} >0$
    depending only on $B$ and $\sup_{t\ge 0} \|g+h\|_{L^1 _3}$ such that
      \[ \forall \, t \ge 0, \ \ \
         \|g_t -h_t \|_{L^1 _2} \le  \|g_{\mbox{\scriptsize{{\em in}}}} 
         - h_{\mbox{\scriptsize{{\em in}}}} \|_{L^1 _2} \,
         e^{C_{\mbox{\scriptsize{{\em stab}}}} t}. \]
    \end{prop}

\smallskip\noindent
{\sl Proof of Proposition~\ref{stab}. }
We multiply the equation satisfied by $(g-h)$ by
$\phi(t,v) = \hbox{sgn}(g(t,v)-h(t,v)) \, (1+|v|^2)$.
We use on the one hand the same arguments as
in~\cite[Proposition~3.4]{MMR1}
to treat
    \[ I= \int_{\RR^N} \left[ Q(g,g) - Q(h,h) \right] \, \phi(t,v) \, 
dv, \]
which gives
    \[ I \le  C \, \left( \int_{\RR^N} (g+h) \, (1+|v|^3) \, dv \right)
         \left( \int_{\RR^N} |g-h| \, (1+|v|^2) \, dv \right). \]
On the other hand we use that
  \bean 
  - \int_{\RR^N} \nabla_v \cdot (v \, (g-h)) \, \phi(t,v) \, dv
  &=& -N \int_{\RR^N} |g-h| \, (1+|v|^2) \, dv \\
  && \ \ \ \ + \int_{\RR^N} |g-h| \, \nabla_v \cdot ( v + v |v|^2 ) \, dv \\
  &=& 2 \, \int_{\RR^N} |g-h| \, |v|^2 \, dv.
  \eean
This concludes the proof with
$C_{\mbox{\scriptsize{stab}}} = C \, \sup_{t\ge 0} \|g+h\|_{L^1 _3} +2$. 
\qed

\subsection{Uniform propagation of Lebesgue norms}\label{subsec:Leb}

Let us take a normal restitution coefficient $e \in (0,1)$ 
(the case $e=0$ can be included in dimension $N=3$) and 
$1< p < +\infty$, and let us consider some initial 
datum $g_{\mbox{\scriptsize{in}}} \in L^1 _2 \cap L^p$. 
We compute the time derivative of the $L^p$ norm of the solution 
$g$ to equation~(\ref{eqresc}):
    \[ \frac1p \, {d \over dt} \int_{\R^N} g^p \, dv =  \int_{\R^N} Q^+(g,g) \, 
g^{p-1}\, dv
       -  \int_{\R^N} g^p \, L(g) \, dv
       - \int_{\R^N} g^{p-1} \, \nabla_v (v \, g) \, dv. \]
We use the control (\ref{Lgv}), and
    \[  \int_{\R^N}  \nabla_v \cdot (v g) \, g^{p-1} \, dv = 
N\left(1-\frac1p\right) \, \|g\|_{L^p} ^p. \]
Gathering all these estimates, we deduce
    \[ \frac1p \, {d \over dt} \int_{\R^N} g^p \, dv \le \int_{\R^N} Q^+(g,g) \, 
g^{p-1} \, dv
    - \min\left\{1, N\left(1-\frac1p\right) \right\} \int_{\R^N} g^p (1+|v|) \, 
dv. \]
Concerning the gain term, Theorem~\ref{regQ} yields, for any $\e>0$,
    \[ \int_{\R^N} Q^+ (g,g) f^{p-1} \, dv \le
     C_\e \, \left\|g\right\|_{L^1 _2} ^{1+p\theta} 
\|g\|_{L^p}^{p(1-\theta)}
     + \e \, \left\|g\right\|_{L^1 _2}  \|g (1+|v|) \|_{L^p} ^p. \]
Hence, using the bound $C_E$ on the kinetic energy,
if we fix $\e$ such that
   \[ C_E  ^{p(1-\theta)} \e < \frac12 \,\min\left\{1,N\left(1-\frac1p 
\right)\right\}, \]
we obtain
    \[ \frac{d}{dt} \|g\|_{L^p} ^p \le C_+ \,  \|g\|_{L^p} ^{p(1-\theta)}
    - K_- \, \|g \|_{L^p _{1/p}} ^p \]
for some explicit constants $C_+$, $K_- >0$.
By maximum principle, it shows that the $L^p$ norm of $g$ is uniformly
bounded by
    \[ \sup_{t \ge 0} \|g_t\|_{L^p} \le
            \max \left\{ \left( \frac{C_+}{K_-}
\right)^{\frac{1}{p\theta}}, \| g_{\mbox{\scriptsize{in}}} \|_{L^p} 
\right\}. \]
The proof for weighted $L^p$ norms is exactly similar. 
This shows the part concerned with Lebesgue norms in point~(i) of 
Theorem~\ref{asymp}. 

\subsection{Non-concentration in the rescaled variables and Haff's law}

In this subsection we give a short proof of Haff's law,
even if a stronger pontwise estimate from below on the tail in rescaled 
variables will be proved in the next section. 
Let us take a normal restitution coefficient $e \in (0,1)$ 
(the case $e=0$ can be included in dimension $N=3$). 
Let $f_{\mbox{\scriptsize{in}}} = g_{\mbox{\scriptsize{in}}}$
be an initial datum in $L^1 _2 \cap L^p$ (with $1<p<+\infty$).
Hence according to the previous subsection,
the rescaled solution $g$ to~(\ref{eqresc}) with initial datum 
$g_{\mbox{\scriptsize{in}}}$
satisfies
     \[ \sup_{t \ge 0} \|g_t\|_{L^p} \le C_p \]
for some explicit constant $C_p >0$ depending on the collision kernel
and the mass, kinetic energy and $L^p$ norm of $f_{\mbox{\scriptsize{in}}}$.
By using Cauchy-Schwarz inequality, this non-concentration
estimate implies that for any $r>0$
     \[ \forall \, t \ge 0, \hspace{0.3cm} \int_{|v|\le r} g(t,v) \, dv
        \le C \, r^{\frac{p-1}p N}. \]
Thus there is $r_0>0$ such that
     \[ \forall \, t \ge 0, \hspace{0.3cm} \int_{|v|\le r_0} g(t,v) \, dv
   \le \frac\rho2 \]
and thus
     \bear
     \label{EE>r2} \forall \, t \ge 0, \quad
            \int_{\RR^N} g(t,v) \, |v|^2 \, dv &\ge& \int_{|v|\ge r_0}  
    g(t,v) \, |v|^2 \, dv \\
      \nonumber       &\ge& r_0 ^2 \, \int_{|v|\ge r_0} g(t,v) \, dv \\
      \nonumber      &\ge& r_0 ^2 \, \left(1-\int_{|v|\le r_0} g(t,v) \, 
dv \right) \ge \frac{\rho \, r_0 ^2}2.
     \eear
As a conclusion, gathering (\ref{BorneY2}) and (\ref{EE>r2}), we have 
proved that for some constants $C_0, C_1 \in (0,\infty)$ there holds
  $$
  C_0 \le \EE( g (t,\cdot)) \le C_1,
  $$
and Haff's law (\ref{hafflaw}) follows thanks to (\ref{momentgtof}), 
which proves Theorem~\ref{Haff}.

\begin{rem} The inequality  $\EE(f(t,\cdot)) \le M \, (1+t)^{-2}$ (or 
equivalently $\EE(g(t,\cdot)) \le C_1$) was already known: see for instance  
\cite[equations~(2.5)-(2.6)]{BenedettoCP97} where it is proved 
for a quasi-elastic one-dimensional model with the same evolution 
equation~(\ref{eqdiffEE}) on the kinetic energy, by comparison to a 
differential equation. Indeed the harder part in Haff's law is the first inequality, which 
means that the solution does not cool down faster than the self-similar
profile. As emphasized by the proof above, this is related to the 
impossiblity of asymptotic concentration in the rescaled equation (\ref{eqresc}).
%\smallskip
%
%2. Points~(i) and~(ii) of Theorem~\ref{asymp} also translate into
%corresponding results on the solution
%in the original variables.
\end{rem}

\subsection{Uniform propagation of Sobolev norms}

Let us take a normal restitution coefficient $e \in (0,1)$ 
(the case $e=0$ can be included in dimension $N=3$). 
The study of propagation of regularity and exponential decay
of singularities is based on a Duhamel representation of the solution
we shall introduce. Let us denote
    \[ L(t,v) = %\|b\|_{L^1(\Sp^{N-1})} \,
     \left( \int_{\RR^N} \, g(v_*) \, |v-v_*| \, dv_*\right), \]
and
    \[ S_t g = g(e^{-t} v) \, \exp \left[ -N t - \int_0 ^t 
L(s,e^{-(t-s)}v) \, ds \right] \]
the evolution semi-group associated to %the unbounded operator
    \[ Tg = -  %\|b\|_{L^1(\Sp^{N-1})} \,
            \left( \int_{\RR^N} \, g(v_*) \, |v-v_*| \, dv_*\right) g(v)
            - \nabla_v \cdot (v \, g). \]

Then the solution of~(\ref{eqresc}) represents as
    \[ g_t = S_t g_{\mbox{\scriptsize{in}}} + \int_0 ^t S_{t-s} Q^+(g_s,g_s) \, ds. \]

We give a proposition similar to~\cite[Proposition~5.2]{MV**}:
    \begin{prop}\label{duh}
    There are some constants $\alpha>0$, $\delta>0$, $K>0$ and $k >0$ such that
    for any $s,\eta \ge 0$, we have
      \[ \|S_t g_{\mbox{\scriptsize{{\em in}}}}\|_{H^{s+\alpha} _\eta} \le
           C_{\mbox{\scriptsize{{\em Duh}}}} \, e^{-Kt} \,
            \|g_{\mbox{\scriptsize{{\em in}}}}\|_{H^{s+\alpha} _{\eta+\delta}}
           \,  \sup_{0 \le \bar t \le t} \|g(\bar t, \cdot)\|_{H^s _{\eta+\delta}} ^{s+k} \]
      \[ \left\| \int_0 ^t S_{t-s} Q^+(g_s,g_s) \, ds \right\|_{H^{s+\alpha} _\eta} 
      \le C_{\mbox{\scriptsize{{\em Duh}}}} \,
           \sup_{0 \le \bar t \le t} \|g(\bar t, \cdot)\|_{H^s _{\eta+\delta}} ^{s+k}. \]
    \end{prop}

\smallskip\noindent
{\sl Proof of Proposition~\ref{duh}.}
The proof is exactly similar to~\cite[Proof of Proposition~5.2]{MV**}.
Indeed
the semi-group in~\cite[Proof of Proposition~5.2]{MV**} is
  \[ \bar{S}_t g = g(v) \, \exp \left[ - \int_0 ^t L(s,v) \, ds \right] \]
and thus the estimates on the Sobolev norm in $v$ can only improve for
$S_t$ according to $\bar{S}_t$. The main tool of~\cite[Proof of
Proposition~5.2]{MV**}, {\it i.e.}, the Bouchut-Desvillettes-Lu regularity 
result on $Q^+$, has been
proved in our case in Theorem~\ref{regQBDL}. \qed

\medskip

Now results follow as in~\cite{MV**}:

   \begin{theo}\label{prop:reg}
   Let $0 \le g_{\mbox{\scriptsize{{\em in}}}} \in L^1_2$ be an initial datum and let $g$ be the unique
   solution of~(\ref{eqresc}) in $C(\R_+,L^1 _2) \cap L^1(\R_+,L^1 _3)$ associated with
   $g_{\mbox{\scriptsize{{\em in}}}}$.
   Then for all $s>0$ and $\eta \ge 1$,
   there exists $w(s)>0$ (explicitly $w(s)=\delta \lceil s/\alpha \rceil$,
   where $\alpha$ is defined in Proposition~\ref{duh}) such that
    \[ g_{\mbox{\scriptsize{{\em in}}}} \in H^s_{\eta+w} \Longrightarrow
    \sup_{t\geq 0} \|g(t,\cdot)\|_{H^s _\eta} < +\infty\]
   with uniform bounds.
   \end{theo}

\smallskip\noindent
{\sl Proof of Theorem~\ref{prop:reg}.}
Let $n \in \N$ be such that $n \alpha \geq s$ ($n=\lceil
s/\alpha \rceil$). Let $w(s)=\delta \lceil s/\alpha \rceil$. The
proof is made by an induction comprising $n$ steps,
proving successively that $g$ is uniformly bounded in
$H^{i \alpha} _{\eta + \frac{n -i}{n} w}$
for $i= 0,1,\dots,n$. %The above-mentioned argument is used in each step.

Let us write the induction. The initialisation for $i=0$, {\em i.e.}, $g$
uniformly bounded in $L^2 _{\eta + w}$ is proved by the previous
study of uniform propagation of weighted $L^p$ norms in 
Subsection~\ref{subsec:Leb}.
Now let $0 < i \le n$ and suppose the induction assumption 
to be satisfied for all $0 \le j < i$.
Then proposition~\ref{duh} implies
   \[
   \|S_t g_{\mbox{\scriptsize{in}}} \|_{H^{i \alpha} _{\eta + \frac{n -i}{n} w}}
   \le C_{\mbox{\scriptsize{Duh}}} \, e^{-Kt}
   \left\|g_{\mbox{\scriptsize{in}}}(\cdot) \right\|_{H^{i \alpha} _{\eta + \frac{n - i}{n} w
+ \delta}} \, \sup_{ 0 \le t_0 \le t}
   \left\| g(t_0,\cdot) \right\|_{H^{(i-1)\alpha} _{\eta + \frac{n -i}{n} w + \delta}} ^{i  \alpha +k},
   \]
and
   \[
   \left \|\int_0 ^t S_{t-s} Q^+(g_s,g_s) \, ds \right\|_{H^{i \alpha} _{\eta + \frac{n -i}{n}w}}
   \le  C_{\mbox{\scriptsize{Duh}}} \, \sup_{0 \le t_0 \le t}
   \left\| g(t_0,\cdot) \right\|_{H^{(i-1)\alpha}
   _{\eta + \frac{n -i}{n}w + \delta}} ^{i \alpha +k} .
   \]
Moreover as $i \ge 1$,
   \[ \eta + \frac{n-i}{n}w + \delta \le \eta + \frac{n -(i-1)}{n} w. \]
Thus, using the induction assumption for $i-1$, $g$ is uniformly
bounded in $H^{i \alpha} _{\eta
+\frac{n-i}{n}w}$, which concludes the proof. \qed

\subsection{Exponential decay of singularities}

Let us take a normal restitution coefficient $e \in (0,1)$ 
(the case $e=0$ can be included in dimension $N=3$). 
In this part we shall follow a similar strategy as in~\cite{MV**} in order to 
show that singularities decrease exponentially fast along the flow in 
rescaled variables. Namely we prove the

   \begin{theo}\label{MM:theo:dec}
   Let $0 \le g_{\mbox{\scriptsize{{\em in}}}} \in L^1 _2 \cap L^p$, $p \in (1,\infty)$,  
   and let $g$ be the unique solution of~(\ref{eqresc}) in 
   $C(\R_+,L^1 _2) \cap L^1 _{\mbox{{\scriptsize {\em loc}}}} (\R_+,L^1 _3)$ associated with
   $g_{\mbox{\scriptsize{{\em in}}}}$. Let  $s\geq 0$, $q\geq 0$ be 
   arbitrarily large. Then
   $g$ can be written $g^S+g^R$ in such a way that
     \[
    \left\{
     \begin{array}{l}
     \displaystyle \sup _{t \ge 0}
     \left\|g_t ^S \right\|_{H^s_q \cap L^1 _2} < +\infty , \quad g^S \ge 0\vspace{0.3cm} \\
     \displaystyle \exists  \, \lambda>0; \ \ \left\|g_t ^R \right\|_{L^1 _2}   = O\left(e^{-\lambda t}\right) .
     \end{array}
    \right.
    \]
   All the constants in this theorem can be computed in terms of
   the collision kernel, the mass and kinetic energy and $L^2$ norm 
   of $g_{\mbox{\scriptsize{{\em in}}}}$.
   \end{theo}

\smallskip\noindent
{\sl Proof of Theorem~\ref{MM:theo:dec}.}
%It is straightforward that the statement~\cite[Theorem~5.4]{MV**} implies 
%the statement of Theorem~\ref{MM:theo:dec}. 
Assume first that $0 \le g_{\mbox{\scriptsize{in}}} \in L^1 _2 \cap L^p$, $p \in [2,\infty)$. 
Then $g_{\mbox{\scriptsize{in}}} \in L^2$ and the proof of Theorem~\ref{MM:theo:dec} 
is exactly similar to~\cite[Proof of Theorem~5.5]{MV**} since
the only tools of the proof are the stability result, the
estimate on the Duhamel representation and the uniform propagation
of Sobolev norms, %and the non-concentration estimate 
%on the iterated gain term and the deduced Abrahamsson's decomposition, 
which have been proved respectively in Proposition~\ref{stab},
Proposition~\ref{duh} and Proposition~\ref{prop:reg}. %Proposition~\ref{itarteQ+} 
%and the discussion after this proposition.
The propagation and appearance of moments in $L^1$ (used in this proof) 
were proved in Proposition \ref{MM:Y3uniform}. 
Moreover as was already pointed out in~\cite[Section~7, Remark~3]{MV**}, 
it is possible with the same arguments 
to relax the assumptions on the initial datum to 
$0 \le g_{\mbox{\scriptsize{in}}} \in L^1 _2 \cap L^p$ for any 
$p \in (1,2]$ by using the gain of integrability of the gain part 
of the collision operator. \qed

%This theorem proves
%the regularity part of point~(ii) in
%Theorem~\ref{selfsim} and
\smallskip
Point~(i) of Theorem~\ref{asymp} is deduced from this theorem.

\begin{rem}
We do not know how to carry 
the argument in~\cite[Theorem~7.2]{MV**} in order to reduce the 
assumptions to only $0 \le g_{\mbox{\scriptsize{in}}} \in L^1 _2$. The 
estimates of Abrahamsson on the iterated gain term can be easily 
extended to the inelastic framework, but the decomposition 
of Abrahamsson (quoted in the elastic case in~\cite[Lemma~7.1]{MV**}) 
between a part with finite $L^p$ norm for some 
$p\in (1,3)$ and a part which decreases exponentially fast requires 
a lower bound on the energy. Here for the rescaled inelastic problem 
we deduce this lower bound from the propagation of $L^p$ bounds, which therefore seem 
compulsory in our method. 
\end{rem}

\begin{rem}
A suggested by this study, the self-similar variables
are not only useful for proving the existence of self-similar profiles, 
but it seems that they also provide the good framework for studying precisely the 
regularity of the solution. For instance, coming back to the original variables, 
Theorem~\ref{MM:theo:dec} shows the algebraic decay of singularities for 
the solutions of (\ref{eqB1}). 
\end{rem}

%%%%%%%%%%%%%%%%%%%%%%%%%%%%%%%%%%%%%%%%%%%%%%

\section{Self-similar solutions and tail behavior}\label{sec:tail}
\setcounter{equation}{0}
\setcounter{theo}{0}

%%%%%%%%%%%%%%%%%%%%%%%%%%%%%%%%%%%%%%%%%%%%%%

In this section we achieve the proofs of Theorem~\ref{selfsim} and
Theorem~\ref{asymp} by
showing the existence of self-similar solutions, and obtaining estimates
on their tail and the tail of generic solutions. We consider 
a normal restitution coefficient $e \in (0,1)$ (and as before 
the case $e=0$ can be included in dimension $N=3$). 

\subsection{Existence of self-similar solutions}

The starting point is the following result, see for instance
\cite[Theorem~5.2]{GPV**} or \cite{Balabane,EMRR}.

     \begin{theo}\label{GPV}
     Let ${\cal Y}$ be a Banach space and $(S_t)_{t \ge 0}$ be a continuous
     semi-group on ${\cal Y}$. Assume that there exists ${\cal K}$ a nonempty convex and
     weakly (sequentially) compact subset of ${\cal Y}$ which is invariant
     under the action of $S_t$  (that is $S_ty  \in {\cal K}$ for any $y \in {\cal K}$
     and $t \ge 0$), and such that $S_t$ is weakly (sequentially) continuous on ${\cal K}$ 
     for any $t>0$. Then there exists $y_0 \in {\cal K}$ which is stationary under the
     action of $S_t$ (that is $S_t y_0=y_0$ for any $t \ge 0$).
     \end{theo}

\smallskip\noindent
{\sl Proof of Theorem~\ref{selfsim} (existence part).}
The existence of self-similar solutions follows from the
application of this result to the evolution semi-group
of~(\ref{eqresc}). The continuity properties of the
semi-group are proved by the study of the Cauchy problem, recalled in
Section~\ref{sec:self-sim}. On the Banach space ${\cal Y} = L^1 _2$, 
thanks to the uniform bounds on the $L^1 _3$ and $L^p$ norms, 
the nonempty convex subset of ${\cal Y}$
     \[ {\cal K} = \left\{ 0 \le f \in {\cal Y}, \hspace{0.3cm} 
         \int_{\RR^N} f \, \left( \begin{array}{ll} 1 \\ v \end{array} \right) \, dv =  
         \left( \begin{array}{ll} \rho \\ 0 \end{array} \right) \quad \mbox{ and } 
         \quad \|f\|_{L^1 _3} + \|f\|_{L^p} \le M \right\} \]
is stable by the semi-group provided $M$ is big enough. This set is
weakly compact in ${\cal Y}$ by Dunford-Pettis Theorem, and the continuity 
of $S_t$ for all $t \ge 0$ on ${\cal K}$ follows from Proposition~\ref{stab}.  
This shows that there exists a nonnegative stationary solution to~(\ref{eqresc}) in $L^1 _3 \cap L^p$ 
for any given mass, that is a self-similar solution for the original problem~(\ref{eqB1}).

Then one can apply Theorem~\ref{MM:theo:dec}, 
which proves that the stationary solution of~(\ref{eqresc}) obtained above 
belongs to $C^\infty$  (in fact it proves that
it belongs to the Schwartz space of $C^\infty$ functions decreasing faster than 
any polynomial at infinity). 
Moreover, since the property of being radially symmetric is stable
along the flow of (\ref{eqresc}), this sationary solution can be shown to exist 
within the set of radially symmetric functions by the same arguments.
\qed

\subsection{Tail of the self-similar profiles}

In this subsection we prove pointwise bounds on the tail behavior of
the self-similar solutions. The starting point is the following result
extracted from~\cite[Theorem~1]{BGP**}; notice that it is also a consequence of 
the construction of invariant sets $\Cc_x$ for $z_p$ with $a =2$,   as 
defined in (\ref{MM:superS}).

     \begin{theo}[Bobylev-Gamba-Panferov] \label{BGPtail}
     Let $G$ be a steady state of~(\ref{eqresc}) with finite moments of all orders.
     Then $G$ has exponential tail of order $1$, that is
       \[ r^* = \sup \left\{ r \ge 0, \quad \int_{\R^N} G(v) \, \exp(r|v|)
          \, dv < +\infty \right\} \]
     belongs to $(0,+\infty)$.
     \end{theo}
Note that if one define more generally (for $s>0$)
     \[  r^* _s = \sup \left\{ r \ge 0, \quad \int_{\R^N} G(v) \,
         \exp(r|v|^s) \, dv < +\infty \right\}, \]
a simple consequence of this result is that $r_s^*=+\infty$ for any $s<1$, 
and $r_s ^* = 0$ for any $s > 1$.

First let us prove the pointwise bound from above on the steady state.
Since the evolution equation (\ref{eqresc}) makes all the moments appear
(see Proposition \ref{MM:Y3uniform}),
we assume that $G$ has finite moments of all orders.
Moreover, as discussed above, we can also assume
that $G$ is smooth and radially symmetric.
We denote $r=|v|$. We thus have the
    \begin{prop}\label{upbdprof}
    Let $G \in C^1$ be a radially symmetric nonnegative steady state of~(\ref{eqresc})
    with finite moments of all order.
    Then there exists $A_1,A_2 > 0$ such that
      \[
      \forall \, v \in \R^N, \ \ \ G(v) \le A_1 \, e^{-A_2 \, |v|}.
      \]
    \end{prop}

\smallskip\noindent
{\sl Proof of Proposition~\ref{upbdprof}}.
The differential equation satisfied by $G=G(r)$ writes
     \[ Q(G,G) - N \, G - r \, G' =0. \]
Since $G$ is smooth and integrable, it goes to $0$ at infinity.
By integrating this equation between $r=R$ and $r=+\infty$, we obtain
     \[ G(R) = N \, \int_R ^{+\infty} \frac{G(r)}{r} \, dr
         -  \int_R ^{+\infty} \frac{Q(G,G)}{r} \, dr. \]
One deduces the following upper bound
     \[ G(R) \le N \, \int_R ^{+\infty} \frac{G(r)}{r} \, dr +
            \int_R ^{+\infty} \frac{Q^-(G,G)}{r} \, dr. \]
Since $Q^-(G,G)=G \, (G * \Phi)$, we have
     \[ Q^-(G,G)(v) \le C \, (1+|v|) \, G. \]
Hence, taking $R \ge 1$ leads to
     \[ G(R) \le C \, \int_R ^{+\infty} G(r) \, r^{N-1} \, dr. \]
Finally, since we have by Theorem~\ref{BGPtail}
     \[ \int_0 ^{+\infty} G(r) \, \exp(A_2 \, r) \, r^{N-1} \, dr \le A_0
< +\infty \]
for some constants $A_0,A_2 >0$, we deduce that
     \[ G(R) \le C \, \int_R ^{+\infty} G(r) \, r^{N-1} \, dr \le C \,
A_0  \, \exp(-A_2 \, R)
             = A_1 \, \exp(-A_2 \, R). \]
This concludes the proof. \qed
\smallskip

For the pointwise lower bound, we give here a proof based on a maximum 
principle argument, inspired from the works \cite{GPV**,GPVbis**}. 
We shall in the next subsection give a more general result for 
generic solutions of~(\ref{eqresc}), based on the spreading 
effect of the gain term and the dispersion (or transport) 
effect of the evolution semi-group of  (\ref{eqresc}) (due to the anti-drift term)
in the spirit of \cite{Ca32,PW97}.

     \begin{prop} \label{lobdprof}
     Let $G \in C^1$ be a nonnegative steady state of~(\ref{eqresc}) 
     with finite moments of orders $0$ and $2$ and which is not identically equal to $0$.
     Then there exists $a_1,a_2 > 0$ such that
       \[
       \forall \, v \in \R^N, \ \ \ G(v) \ge a_1 \, e^{-a_2 \, |v|}.
       \]
     \end{prop}

We first start with a lemma.
     \begin{lem}\label{lemlbdprof} For any $r_0,a_1, \rho_0,  \rho_1 > 0$, 
     there exists $a_2 > 0$ such that
     the function $h(v) := a_1 \, \exp(-a_2 \, |v|)$ satisfies
       \beqn \label{ingh}
       \forall \, v, \,\, |v| \ge r_0, \quad Q^-(g,h) + \nabla_v \cdot (v \, h) \le 0
       \eeqn
     for any function $g$ such that
       \[ \int_{\R^N} g (v) \, dv = \rho_0, \quad \int_{\R^N} g (v) \, 
          |v| \, dv = \rho_1. \]
     \end{lem}

\smallskip\noindent
{\sl Proof of Lemma~\ref{lemlbdprof}}. On the one hand, it is
straightforward that
    \[ Q^-(g,h) := (g * \Phi) \, h \le (\rho_1 + \rho_0 \, |v|) \, h. \]
On the other hand, simple computations show that
    \[ \nabla_v \cdot (v \, h) = (N - a_2 \, |v|) \, h. \]
Gathering these two inequalities there holds
    \[ \forall \, v, \,\, |v| \ge r_0, \quad
       Q^-(g,h) + \nabla_v \cdot (v \, h) \le (\rho_1 + N + \rho_0 \, |v| -  a_2
       \, |v|) \, h \le 0 \]
for $a_2$ large enough. \qed

\smallskip\noindent
{\sl Proof of Proposition~\ref{lobdprof}}.
Since $G \in C^1$ and it is radially symmetric, there holds $G'(0) = 0$.
As a consequence, the equation satisfied by $G$ reads in $v = 0$
    \[ Q(G,G) (0) - N \, G(0) = 0 \]
and then
    \[ G(0) = {Q^+(G,G)(0) \over N}> 0 \]
since $G$ is not zero everywhere.
By continuity, $G(v) >  2 \, a_1$ on $B(0,r_0)$ for some $a_1,r_0 > 0$.

Let us define
    \[ \rho_0 := \int_{\R^N} G (v) \, dv, \quad \rho_1 :\int_{\R^N} G (v) \, |v| \, dv, \]
and $a_2>0$ associated to $r_0, a_1, \rho_0, \rho_1$ 
by Lemma~\ref{lemlbdprof}. On the one hand
$h (v) :=  a_1 \, \exp(-a_2 \, |v|)$ satisfies~(\ref{ingh}) for $g = G$ and,
on the other hand, $G$ satisfies
    \beqn \label{ingG}
    \forall \, v \in \R^N, \quad Q^-(G,G) + \nabla_v (v \, G) = Q^+(G,G) \ge  0.
    \eeqn
Introducing the auxiliary function $W := G-h$, we deduce
from~(\ref{ingh}) and~(\ref{ingG})
    \[ \forall \, v, \,\,|v| \ge r_0,  \ \ \ (G * \Phi) \, W + \nabla_v (v \, W) \ge 0 \]
and $W(r_0) = G(r_0) - h(r_0) \ge G(r_0)/2 > 0$. By the Gronwall Lemma
(using that all the functions involved in this inequality are radially
symmetric), we get $W(v) \ge 0$  for any $v$, $|v|\ge r_0$, which
concludes the proof. \qed

\subsection{Positivity of the rescaled solution}

\medskip
We start with three technical lemmas.

    \begin{lem}\label{lower1} 
    Let $g_0$ satisfies for $p \in (1,\infty)$
     \beqn \label{condg0}
     \int_{\R^N} g_0 \, dv = 1, \quad
     \int_{\R^N} g_0 \, |v|^2 \, dv \le C_1, \quad
     \int_{\R^N} g_0^p\, dv \le C_2.
     \eeqn
There exist $R > r > 0$ and $\eta > 0$ depending only on $C_1, C_2$, and
$(v_i)_{i=1,\dots,4}$ such that $|v_i|\le R$, $i=1,\dots,4$, and 
$|v_i - v_j|\ge 3r$ for $1\le i \not= j \le 3$, and
\beqn \label{intBrg}
\int_{B(v_i,r)} g_0 (v) \, dv \ge \eta \quad\hbox{ for }\quad i=1,\,2,\,3,
\eeqn
\beqn \label{sphereplan}
\forall \, w_i \in B(v_i,r), \quad E^e_{w_3,w_4} \cap S^e_{w_1,w_2}
\hbox{ is a sphere of radius larger than } r,
\eeqn
where $E^e_{v,'v}$ stands for the plane defined in 
Proposition~\ref{carlrep}
and $S^e_{v,v_*}$ stands for the sphere of all possibles post-collisional 
velocity $v'$ defined by (\ref{vprime}).
\end{lem}

\smallskip\noindent
{\sl Proof of Theorem~\ref{lower1}}.
Let $C_R$ denotes the hypercube $[-R,R]^N$ centered at $v = 0$ with length $2R > 0$.
Thanks to the mass condition and the energy bound in (\ref{condg0}), for
$R$ large enough, there holds
\beqn \label{massCR}
\int_{C_R}g_0 \, dv \ge {1 \over 2}.
\eeqn
Then we define $(K_i)_{i = 1,\dots, I}$ the familly of $I = (2 \,
R/r)^N$ hypercubes of length $r  > 0$ (with $R/r \in \N$), included  in $C_R$
and such that the union of $K_i$ is almost equal to $C_R$. For any
given $\lambda >0$ to be later fixed, we may find $r > 0$ such that
\beqn\label{massKk+r}
\quad \int_{K_i + B(0,\lambda r)} \!\!\! \!\!\!  \!\!\! g_0 \, \, dv
   \le |K_i + B(0,\lambda  r) |^{1/p'} \left( \int_{K_i + B(0,\lambda
    r)} \!\!\!\!\!\! \!\!\! g_0^p  \, \, dv \right)^{1/p} \le C \,
[(\lambda+1)  r] ^{N/p'} \le {1/4}
\eeqn
for any $i = 1,\dots, I$. Hence we can choose $K_{i_0}$ such that the 
mass of $g_0$ in $K_i$ is maximal for $i=i_0$. Because of (\ref{massCR}) there
holds
   \beqn\label{massKiKj}
   \int_{K_{i_0}} \, g_0 \, dv \ge {1/4} \, (2 \, R/r)^{-N}. 
%\quad\hbox{ for
%}\quad k=i.
   \eeqn
Gathering (\ref{massCR}) and (\ref{massKk+r}) we may find $K_{j_0} \subset
C_R$ such that dist$(K_{i_0},K_{j_0}) > \lambda r$ and (\ref{massKiKj}) also
holds for $i = j_0$.

\smallskip
Next, we fix $\lambda := 200 \, \beta$. We define 
$v_1$ (respectively $v_2$) as the center of the hypercube 
$K_{i_0}$ (respectively $K_{j_0}$), and 
$v_3 = (v_1+v_2)/2$ and $v_4 = v_2$. 
Then we have 
  \[ \Omega(v_3,v_4) = v_1 + {\beta^{-1} \over 2} \, (v_2 - v_1) \in [v_1,v_2], \] 
which implies 
  \[ |\Omega - v_1| = {\beta^{-1} \over 2} \, |v_2 - v_1| \ge {\beta^{-1} \over 2}
     \, (\lambda \, r) \ge 100 \, r. \] 
Thus $E^e_{v_3,v_4} \cap S^e(v_1,v_2) $ is a $(N-2)$-dimensional sphere of radius
larger than $100 \, r$ (because $B(\Omega,100r)$ is included in the convex hull of 
$S^e(v_1,v_2)$), and (\ref{sphereplan}) follows straightforwardly. \qed

\begin{lem}\label{lower2} 
Let us fix $R > r > 0$ and $\eta > 0$. Then there 
exists $\delta_0 >0$, $\eta_0 > 0$, $\xi_0 \in (0,1)$ (depending on $R > r > 0$, 
$\eta  > 0$ and $B$) such that, for any functions $f$, $h$, $\ell$
satisfying (\ref{intBrg})-(\ref{sphereplan}) for some velocities 
$(v_i)_{i=1,\dots,4}$ such that $|v_i|\le R$, $i=1,\dots,4$ and 
$|v_i - v_j|\ge 3r$, $1 \le i \not= j \le 3$, 
and for any $\xi \in (\xi_0,1)$, there holds
   \[ %\label{QpQpinf}
   Q^+ ( f,  Q^+_\xi ( h, \ell) )  \ge \eta_0 \, {\bf 1}_{B(v_3,\delta_0)} ,
   \]
where we define here and below $Q^+_\xi(\cdot,\cdot) (v) = Q^+(\cdot,\cdot) (\xi \, v)$.
\end{lem}

\smallskip\noindent
{\sl Proof of Theorem~\ref{lower2}}.  We first establish a convenient  
formula to handle representations of the iterated 
gain term. For any $f$, $h$ and $\ell$ and any $v \in \RR^N$
there holds (setting $'v = w$ and $'v_* = w_*$)
\bean
Q^+(f,Q^+_\xi(h,\ell)) (v)
= C'_b \int_{\RR^N}  {f(w) \over |v-w|} \left\{ \int_{E^e _{v,w}}
Q^+_\xi(h,\ell)(w_*)  \, dw_* \right\}\, dw.
\eean
From the following identity
$$
Q^+_\xi(h,\ell)(w_*) =  \int_{\RR^N} \! \int_{\RR^N}  h(w_1) \, \ell(w_2)
\, Q^+_\xi(\delta_1,\delta_2) (w_*) \, dw_1 dw_2
$$
where $\delta_i$ stands for the Dirac measure at $w_j$, the
term between brackets, that we denote by $A$, write
$$
A(v,w) =  \int_{\RR^N \times \RR^N}  h(w_1) \, \ell(w_2) \,
\left\{ \lim_{\eps\to0} {1\over 2 \, \eps} \int_{\RR^N}
Q^+(\delta_1,\delta_2) (\xi \, w_*)  \, \Xi_\eps(w_*) \, dw_* \right\} \,
dw_1 \, dw_2
$$
where $\Xi_\eps$ denotes the indicator function of the set 
$\{w_* \, ; \ \hbox{dist}(w_*,E_{v,w})< \eps \}$.
Denoting now by $D_\eps$ the integral  just after the limit sign in the
term between brackets, and using the weak formulation
(\ref{Qplusweak}), there holds
\bean
D_\eps &=& {\xi^{-1} \over 2 \, \eps} \, \int_{\RR^N \times \RR^N \times \Sp^{N-1}} 
\delta_1(z) \,  \delta_1(z_*) \, |z-z_*|\, b(\sigma
\cdot \hat z) \, \Xi_\eps(\xi^{-1} \, z') \,  d\sigma \, dz \, dz_* \\
&=& |w_1 - w_2| \, \xi^{-1} \, C_b \int_{\Sp^{N-1}} {\Xi_\eps(z'\,
\xi^{-1} ) \over 2\eps}\,  d\sigma,
\eean
where in these integrales $z'$ is defined from $(z,z_*,\sigma)$ and
next from $(w_1,w_2,\sigma)$ thanks to formula (\ref{vprime}). We
define $\xi_0 = (1+r/R)^{-1}$ in such a way that $|\xi^{-1} \, z'-z'|
\le r$ for any $z' \in B(0,R)$ and $\xi \in (\xi_0,1)$. Taking $v \in B(v_3,r)$, $w \in
B(v_4,r)$, $w_1 \in B(v_1,r)$, $w_2 \in B(v_2,r)$, we have thanks to
(\ref{sphereplan}) and $|w_2-w_1|\ge r$:
\[ %\label{infDeps}
D_0 (v,w,w_1,w_2) := \lim_{\eps \to 0} D_\eps \ge r \,  \xi^{-1}_0 \,
C_b \, C \, r^{N-2}.
\]
As a consequence, for any $v \in B(v_3,r)$,
\bean
Q^+(f,Q^+_\xi(h,\ell)) (v) &\ge& Q^+(f \, {\bf 1}_{B(v_4,r)} ,Q^+(h \,
{\bf 1}_{B(v_1,r)},\ell \, {\bf 1}_{B(v_2,r)})) (v)  \\
&\ge& C'_b \int_{B(v_1,r)} \! \int_{B(v_2,r)} \! \int_{B(v_4,r)}
{f(w) \over |v-w|} \, h(w_1) \, \ell(w_2) \,
D_0 \,   dw_1 \, dw_2 \, dw \\
&\ge& C'_b \, \eta^3 \, {1 \over 2 \, R} \,  r \, \xi^{-1}_0 \,   C_b
\, C \, r^{N-2} =: \eta_0.
\eean
This concludes the proof. \qed

\begin{lem}\label{MM:lem:spread} 
For any $\bar{v} \in \R^N$ and $\delta > 0$,  
there exists $\kappa = \kappa(\delta) > 0$ such that
   \begin{equation}\label{MM:eq:spread} 
   \Qq^+ (v) :=  Q^+ \big( {\bf 1}_{B(\bar{v},\, \delta)} , {\bf 1}_{B(\bar{v},\, \delta)}\big)
   \ge \kappa \,
   {\bf 1}_{B\big(\bar{v}, {\sqrt{5} \over 2}\, \delta \big)}.
   \end{equation}
\end{lem}

\medskip
\noindent{\sl Proof of Lemma~\ref{MM:lem:spread}. } The homogeneity 
property (\ref{homogeneiteQ}) of $Q^+$ and the invariance by 
translation allow to reduce the proof of~(\ref{MM:eq:spread}) to the case 
$\bar v = 0$ and $\delta = 1$. The invariance by rotations implies that 
$\Qq^+ $ is radially symmetric and the homogeneity property again 
allows to conclude that the support of $\Qq^+$ is a ball $B'$. More 
precisely, taking a $C^\infty$ radially symmetric  function $\phi$ such 
that $\phi > 0$ on $B = B(0,1)$ and $\phi \le {\bf 1}_B$ on $\R^N$, we 
have $Q^+(\phi,\phi)$ is continuous, $\Qq^+ \ge Q^+(\phi,\phi)$ on 
$\R^N$ and $Q^+(\phi,\phi) > 0$ on the ball $B'$. As a consequence, for 
any ball $B''$ strictly included in $B'$, there exists $\kappa > 0$ 
such that $\Qq^+ \ge \kappa \, {\bf 1}_{B''}$. In order to conclude, we 
just need to estimate the support of $\Qq^+$. 

Let us fix $R \in (0,1)$ and choose $'v ,\! \! ~ 'v_* 
\in B(0,1)$ such that $'v  \perp \! \! ~'v_* $, $|\!~'v| = |\!~'v_*| = R$. Then for 
any $\sigma \in \Sp^{N-1}$, $\sigma \perp \! ~'v- \! ~'v_*$, the function $\Qq^+$ 
is positive at the post-collisional associated velocity $v$ defined by
$$
v = {'v+\! \!~'v_* \over 2} + {1-e \over 4} \, (\!~'v- \! \!~'v_*) +  {1+e \over 4} \, 
|\!~'v-\! \!~'v_*| \, \sigma.
$$
Remarking that $|\!~'v+\!\!~'v_*|^2 = |\!~'v-\! \!~'v_*|^2 = 2 \, R^2$, $(\!~'v-\!\!~'v_*) \cdot 
(\!~'v+\! \!~'v_*) = 0$ and $(\!~'v+\! \!~'v_*) \cdot \sigma = \sqrt{2} \, R$, 
we easily compute
$$
|v|^2 = R^2 \, \left[ 1 + \left( {1+e \over 2} \right)^2 \right] > {5 
\over 4}\, R^2,
$$
and the radius of $B'$ is strictly larger than $\sqrt{5}/2$.
\qed

\medskip

  \begin{theo}\label{theopositivity} Let $g_{in}$ satisfy the hypothesis 
  of Theorem~\ref{asymp} and let $g$ be the solution to the rescaled equation 
  (\ref{eqresc}) associated to  the initial datum $g_{in}$. Then for any $t_*>0$, 
  $g(t,\cdot) > 0$ a.e. on $\RR^N$ for any $t \ge t_*$, 
  and there exists $a_1,a_2, c > 0$ such that
    \[
    \forall \, t \ge t_*, \quad g(t,v) \ge 
    a_1 \, e^{-a_2 \, |v|} \, {\bf 1}_{|v| \le c \, e^{t-t_*}} \ \ \mbox{ for a.e. } v \in \R^N.
    \]
  \end{theo}

\medskip
\noindent{\sl Proof of Theorem~\ref{theopositivity}. }  We split the 
proof into four steps.

\smallskip\noindent
{\sl Step 1.} The starting point is the evolution equation satisfied by $g$ 
written in the form
$$
     \partial_t g + v \cdot \nabla_v g + (N + |v|)  \, g = Q^+(g,g) + (|v| - L(g))  \, g.
$$
Let us introduce the semi-group $S_t$ associated to the operator 
$v \cdot \nabla_v + \lambda(v)$, where $\lambda (v) := N+|v|$. Thanks to the
Duhamel formula and (\ref{Lgv}), we have
\beqn\label{DuhamelF}
g(t,\cdot) \ge S_t \, g(0,\cdot) + \int_0^t S_{t-s} Q^+(g(s,\cdot),g(s,\cdot)) \,ds,
\eeqn
where the semi-group $S_t$ is defined by
\[ %\label{semigrouplambda}
(S_t \, h) (v) =  h(v \, e^{-t}) \, \exp \left( - \int_0^t \lambda(v\,
e^{-s}) \, ds \right).
\]
Notice that 
\[ %\label{int0tlambda}
\left( - \int_0^t \lambda(v\, e^{-s}) \, ds \right) \ge  - (|v|+ N \, t).
\]

\smallskip\noindent
{\sl Step 2. } Let us fix $t_0 > 0$ and define $\tilde g_0 (t,\cdot) := 
g(t_0 + t, \cdot)$. Using twice the
Duhamel formula (\ref{DuhamelF}), we find
\bear\nonumber
\tilde g_0(t,\cdot) &\ge&  \int_0^t S_{t-s} Q^+\left(\tilde g_0(s,\cdot), 
\int_0^s
S_{s-s'} Q^+(\tilde g_0(s',\cdot),\tilde g_0(s',\cdot)) \, ds' \right) \, ds \\ 
\label{glower2}
&\ge&  \int_0^t  \int_0^s S_{t-s} Q^+\left(S_s \tilde g_0,
S_{s-s'} Q^+(S_{s'} \tilde g_0, S_{s'} \tilde g_0)  \right) \, ds' \, ds.
\eear
We apply now Lemma~\ref{lower1} to $\tilde g_0$ and set $R_0 := 2 R$. 
Since $S_t$ is continuous in
$L^1$, there exists $T_1 > 0$, such that for any $s \in [0,T_1]$, there
holds
\[ %\label{intBrStg}
\int_{B(v_i,r)} S_s(\tilde g_0)(v) \, dv \ge \eta/2 \quad\hbox{ for 
}\quad
i=1,\,2,\, 3,
\]
and $e^{-T_1} > \xi_0$.
For $ v \in B(0,R_0)$ and $t \in [0,T_1]$ we may estimate $S_t h$ from below
in the following way
\[ %\label{semigrouplambda}
   (S_t \, h) (v) \ge \gamma\,  h_{e^{-t}} (v)
\]
for some constant $\gamma = \gamma_{R_0,T_1}$. The bound from below 
(\ref{glower2}) then yields (using Lemma~\ref{lower2})
\bear\nonumber
\tilde g_0(t,\cdot)
&\ge& \gamma^2 \,  \int_0^t  \int_0^s Q^+_{e^{s-t}} \left(S_s
\tilde g_0 ,  Q^+_{e^{s'-s}}(S_{s'} \tilde g_0, S_{s'} \tilde g_0)  \right) \, 
ds' \, ds \\
\nonumber
&\ge& \gamma^2 \,  \int_0^t  \int_0^s \eta_0 \, {\bf 1}_{v \,
e^{s-t} \in B(v_3,r)} \, ds' \, ds.
\eear
We have then proved that there exists $T_1 > 0$ and for any $t_1 \in 
(0,T_1/2]$ there exists $\eta_1 > 0$
such that (for some $\bar{v} \in B(0,R)$)
\[
\forall \, t \in [0,T_1/2], \quad 
\tilde g_1( t, \cdot) := \tilde g_0( t+ t_1, \cdot)  \ge \eta_1 \, {\bf 
1}_{B(\bar v,\delta_1)}.
\]

\smallskip\noindent
{\sl Step 3. } Using again the  Duhamel formula (\ref{DuhamelF}) and 
the preceding step we have
\[ %\label{DuhamelF}
\tilde g_1 (t,\cdot) \ge  \int_0^t S_{t-s} Q^+(\tilde g_1 (s,\cdot),\tilde g_1 
(s,\cdot)) \,  ds.
\]
Thanks to Lemma~\ref{lower2}, on the ball $B(0,R_0)$, there holds
\bean %\label{DuhamelF}
\tilde g_1 (t,\cdot)
&\ge&  \eta_1^2 \int_0^t S_{t-s} Q^+\big({\bf 1}_{B(\bar v,\delta_1)},  
{\bf 1}_{B(\bar v,\delta_1)}\big) \,  ds \\
&\ge&  \eta_1^2 \, \kappa(\delta_1) \, e^{-(R_0 + N \, t)}\int_0^t  
{\bf 1}_{e^{-t} \, v \in B(\bar v,\sqrt{5} \, \delta_1/2)} \,  ds \\
&\ge&  \eta_1^2 \, \kappa(\delta_1) \, e^{-(R_0 + N \, T_1)}\, t \,  
{\bf 1}_{B(\bar v,\sqrt{19} \, \delta_1/4)} 
\eean
on $[0,T_2]$ with $T_2 \in  (0,T_1/2]$ small enough, and then
\[
\tilde g_2( t, \cdot) := \tilde g_1( t+ t_2, \cdot)  \ge \eta_2 \, {\bf 
1}_{B(\bar v,\delta_2)}
\quad\hbox{on}\quad [0,T_2/2]
\]
with $\delta_2 := \sqrt{19} \, \delta_1/4$ and $t_2 \in (0,T_2/2]$ 
arbitrarily small, $\eta_2 > 0$. Repeating the argument we obtain
\[
\tilde g_k( t, \cdot) := g\left( t+ \sum_{i=0}^k t_i, \cdot \right) 
\ge \eta_k \, {\bf 1}_{B(\bar v,\delta_k)}
\,\,\hbox{on}\,\, [0,T_k/2], \quad \hbox{with} \,\, \delta_k := 
(\sqrt{19}/4)^k \, \delta_1
\]
with $k \ge 1$ and some $t_i \in [0,T_i/2]$ arbitrarily small, $\eta_k > 0$.
As a consequence, taking $k$ large enough in such a way that $\delta_k 
R_0$, we get for some explicit constant $\eta_* > 0$ and some   
(arbitrarily small) time $t_* > 0$
\beqn\label{eta*}
\forall \, t_0 \ge 0, \quad g(t_* + t_0,\cdot) \ge \eta_* \, {\bf 1}_{B(0,R)}.
\eeqn

\smallskip\noindent
{\sl Step 4. }  Coming back to the Duhamel formula (\ref{DuhamelF}) 
where we only keep the first term, we have, for any $t_0 \ge 0$,
$$
\forall \, t \ge t_*, \quad g(t_0+t,v) \ge 
\eta_* \, {\bf 1}_{|v| \le R \, e^{t-t_*}} \, \exp ( - |v|- N \, (t-t_*)).
$$
As a consequence, for any $t > t_*$,
\bear\label{gexpav}
g(t,v) &\ge& {\bf 1}_{|v| \le R \, e^{t-t_*}}  \left( \sup_{s \in [0,t-t_*]}  
{\bf 1}_{|v| = R \, e^{s}} \, \exp ( -|v|- N \,  s)\right)  \\ 
\nonumber
&\ge&  {\bf 1}_{|v| \le R \, e^{t-t_*}} \, \left( \sup_{s \in [0,t-t_*]}  {\bf 1}_{|v| = R \, e^{s}} \, 
\right)  \exp ( - |v| - N \, \ln^+(|v|/R) ),
\eear
and we conclude gathering (\ref{eta*}) and (\ref{gexpav}). \qed

\medskip

It is straightforward that Theorem~\ref{theopositivity} implies the 
lower bound in point~(iii) of Theorem~\ref{asymp}.

\section{Perspectives}

As a conclusion, we discuss some possible perspectives arising from our study 
(partial answers to them shall be studied in a forthcoming work~\cite{MMIII}).  
\smallskip

Let us denote
  \bean
  {\cal P} &=& \Big\{ G \in C^\infty, \ G \mbox{ radially symmetric}, \\   
     && \hspace{1cm} \exists \, a_1,a_2,A_1,A_2 >0 \ \ | \ \ 
         a_1 e^{- a_2 |v|} \le G(v) \le A_1 e^{- A_2 |v|} \Big\}.
  \eean

\noindent
{\bf Conjecture 1.} For any mass $\rho >0$, the self-similar profile
$G_\rho$ with mass $\rho$ and momentum $0$ is unique.
\smallskip

If Conjecture~1 is true, the natural conjecture is \\
\noindent
{\bf Conjecture 2. (Strong version)} 
For any initial datum with mass 
$\rho$ and momentum $0$ 
(maybe with some regularity and/or moment assumptions), the associated 
solution satisfies (in rescaled variables)
    \[ g_t \rightarrow_{t \to \infty} G_\rho, \]
where $G_\rho$ is the steady state of~(\ref{eqresc}) with mass $\rho$ and momentum $0$.
%or equivalently
%   \[ f_t \sim_{t \to \infty} F_\rho \]
\smallskip

A relaxed version can be \\
\noindent
{\bf Conjecture 2. (Weak version)} For any initial datum with mass 
$\rho$ and momentum $0$ the associated 
solution satisfies (in rescaled variables)
    \[ g_t = g^S _t + g^R _t \]
with $g^S _t \in {\cal P}$ and $g^R _t \to_{t \to \infty} 0$ in $L^1$.
\medskip

Note that the weak version of Conjecture~2 still makes sense when 
the self-similar profile with mass $\rho$ and momentum $0$ is not unique
and even if there is no 
convergence towards some self-similar profile 
(which could be the case for instance
if the solution in rescaled variables ``oscillates'' 
asymptotically between several self-similar profiles).

\bigskip
\noindent
{\bf{Acknowledgments.}} The second author wishes to thank Giuseppe
Toscani for fruitful discussions. Support by the European network HYKE,
funded
by the EC as contract HPRN-CT-2002-00282, is acknowledged.

\end{document}